\documentclass[12pt]{article}
\usepackage{amsfonts,amsmath,times}

\newtheorem{theorem}{Theorem}[section]

\newtheorem{lemma}[theorem]{Lemma}
\newtheorem{coro}[theorem]{Corollary}
\newtheorem{remark}[theorem]{Remark}

\newcommand{\abs}[1]{\vert#1\vert}
\def\b{{\mathbf b}}

\def\Box{{\vrule height 1.3ex width 1.0ex depth -.2ex} \hskip 1.5 truecm}
\def\C{{\mathbb C}}
\def\eps{\varepsilon}
\renewcommand{\Im}[1]{{\mathrm{Im}\left[#1\right]}}

\def\n{{\mathbf n}}

\newcommand{\norm}[1]{\left\Vert#1\right\Vert}

\def\proof{\medskip\noindent{\sc Proof: \ }\nobreak}
\newcommand{\phivec}{\mbox{\boldmath $\phi$}}
\newcommand{\psivec}{\mbox{\boldmath $\Psi$}}
\def\R{{\mathbb R}}
\renewcommand{\Re}[1]{{\mathrm{Re}\left[#1\right]}}
\def\S{{\mathbb S}}
\def\u{{\mathbf u}}
\def\v{{\mathbf v}}
\def\x{{\mathbf x}}
\def\Y{{\cal Y}}

\title{~\\ON THE CAUCHY PROBLEM \\
FOR A DYNAMICAL EULER'S ELASTICA\\[1cm]}
\author{Almut Burchard, Lawrence E. Thomas\\[1cm]
University of Virginia \\
Department of Mathematics\\
Charlottesville, Virginia 22904\\[1cm]
$\{{\tt burchard,let}\}@{\tt virginia.edu}$\\[1cm]}
\date{December 18, 2001; revised May 28, 2002.\\
To appear in Comunications in PDE.}

\begin{document}

\maketitle

\vfill 
\begin{abstract}
The dynamics for a thin, closed loop inextensible Euler's
elastica moving in three dimensions are considered.  The equations of
motion for the elastica include a wave equation involving fourth order
spatial derivatives and a second order elliptic equation for its
tension. A Hasimoto transformation is used to rewrite the equations 
in convenient coordinates for the space and time derivatives of the
tangent vector.  A feature of this elastica is that it exhibits
time-dependent monodromy. A frame frame parallel-transported 
along the elastica is in general only quasi-periodic, resulting in
time-dependent boundary conditions for the coordinates. This
complication is addressed by a gauge transformation, after
which a contraction mapping argument can be applied. Local
existence and uniqueness of elastica solutions are 
established for initial data in suitable Sobolev spaces.  
\end{abstract}


\newpage

\section{Introduction}
\label{sect:intro}
\setcounter{equation}{0}

\subsection{Statement of the result}

Let $\x(s,t)$ be a smooth closed curve in ${\bf R}^3$,
parametrized by its arc length $0\leq s\leq 2\pi$ and by time~$t$.
Physically, the curve can be thought of as a loop of very thin 
inextensible wire which can move in space.  We assume that the
curve is flexible, with potential energy at time $t$ 
determined by its curvature
$\kappa(s,t)= \abs{\partial_s^2\, \x(s,t)}$,
\begin{equation}\label{poteq}
{\cal V}
\ =\ \frac{1}{2}\int_0^{2\pi}\kappa^2(s,t)\,\frac{ds}{2\pi}\ .
\end{equation}
Its kinetic energy at time $t$ is given by
\begin{equation}\label{kineq}
   {\cal T}= \frac{1}{2}\int_0^{2\pi}
       \abs{ \partial_t\x(s,t)}^2\,\frac{ds}{2\pi}\ .
   \end{equation}
We will refer to a dynamical curve $\x(s,t)$ with this 
potential and kinetic energy as a {\em dynamical elastica}.
This choice of kinetic and potential energy
gives rise to the variational
problem of finding extreme solution curves for the 
{\em Lagrangian} ${\cal L}\ = \ {\cal T}-{\cal V}$
subject to the {\em constraint} that $s$ be in fact arc length.
The constraint can be implemented by adding a suitable Lagrange
multiplier term to the Lagrangian, enforcing the condition that the 
tangent vector $ \u(s,t)= \partial_s \x(s,t)$ be a vector of
fixed modulus one. The variational problem we consider is  then
\begin{equation}\label{lageq}
 \delta \int_0^t\left\{ {\cal L}\bigl(\x(\cdot,t),\partial_t\x(\cdot,t)\bigr) 
         -\frac{1}{2}\int_0^{2\pi} \!\!\lambda(s,t)
     \Bigl(\abs{\partial_s \x(s,t)}^2-1\Bigr) \,\frac{ds}{2\pi}\right\}\, dt 
\ =\ 0\ .
\end{equation}
The variational derivative with respect to $\lambda$ provides the
arclength constraint.
The resulting Euler-Lagrange equations are
\begin{equation}
\label{eq:EL-x}
   \partial_t^2 \x\ =\ -\partial_s^4\x +
\partial_s\left(\lambda\partial_s \x\right) \ ,\qquad
   \abs{\partial_s \x}^2\ =\ 1\ ,
\end{equation}
where, we emphasize, the Lagrange multiplier $\lambda(s,t)$ is a
function of both $s$ and $t$.  
{F}rom this equation, it is apparent that $\lambda(s,t)$ has the
physical significance of being the {\em tension} of the curve at
$s$, at time $t$. Note that $\lambda$ can be negative, signifying
compression.  Eq.~(\ref{eq:EL-x}) is nonlinear in $\x$ due to the
fact that at any time $t$, the tension is itself the solution 
to an inhomogeneous elliptic equation involving the tangent 
vector $\u(s,t)=\partial_s\x(s,t)$ and its derivatives 
(see Eq.~(\ref{eq:lambda}) below). We note that
omitting or relaxing the arclength  constraint would require
computing the curvature as 
$\kappa = \partial_s \x \times \partial_s^2 \x / |\partial_s \x|^3$,
leading to a rather forbidding quasilinear Euler-Lagrange equation.

Differentiating Eq.~(\ref{eq:EL-x}) with respect to $s$, one
obtains  for the tangent vector 
the equation
\begin{equation}
\label{eq:EL-u}
\partial_t^2\u\  =\ \partial_s^2\left(-\partial_s^2 + \lambda\right)\u\ ,
\qquad   \u \in \S^2\ .
\end{equation} 
We study here the Cauchy problem for the wave
equation Eq.~(\ref{eq:EL-u}), assuming periodic boundary conditions
\begin{equation}
\u(s\!+\!2\pi, t)\ =\ \u(s,t)\ ,
\label{eq:BC-u}
\end{equation}
and initial conditions
\begin{equation}
\label{eq:initial-u}
\u(s,0)=\u_0(s)\ ,\quad 
\partial_t\u(s,0)=\u_1(s)\ ,
\end{equation}
where $\u_0$ and and $\u_1$ are given $2\pi$-periodic functions with  
\begin{equation}
\label{eq:compatible-u}
    \abs{\u_0(s)}=1\ , \quad \u_0(s)\cdot \u_1(s)=0\
\end{equation}
for all $s$, that is, $\u_0(s)$ lies in $\S^2$, and $\u_1(s)$ is
a tangent vector to the sphere at the point $\u_0(s)$.
We are interested in {\em weak solutions}
of Eq.~(\ref{eq:EL-u}), defined by the property that 
\begin{eqnarray}
\frac{d}{dt} \int_0^{2\pi} 
      \!\phivec(s)\cdot \partial_t\u(s,t)\, \frac{ds}{2\pi}
&=&   \int_0^{2\pi} \!\partial_s^2\phivec(s)\cdot 
        \bigl(- \partial_s^2 \u(s,t) + \lambda(s,t) \u(s,t)\bigr)\, 
\frac{ds}{2\pi}\quad
\label{eq:u-weak}
\end{eqnarray}
for any smooth $2\pi$-periodic function $\phivec$.
We require the map $t\mapsto (\u,\partial_t\u)$ to be strongly
continuous in a suitable Sobolev space, and interpret
the time derivative on the left hand side in the sense of distributions
(see Section~\ref{sect:estimates} for a definition of the Sobolev spaces). 
The following is our main result.

\begin{theorem} Let $(\u_0,\u_1)$ be a pair of $2\pi$-periodic 
functions in a Sobolev space $H^{r+2}\times H^{r}$ with
values in $\R^3$ which satisfy the constraint
in Eq.~(\ref{eq:compatible-u}).  For $r\ge 1/2$, there exists a time $T>0$,
which depends on $\norm{\u_0}_{H^{r+2}}$
and on $\norm {\u_1}_{H^r}$, such that the initial-value problem
given by Eqs.~(\ref{eq:EL-u})-(\ref{eq:initial-u})
has a strongly continuous solution $\u$ on $[0,T]$ 
with $\u(\cdot, t)\in H^{r+2}$, $\partial_t \u(\cdot,t)\in H^{r}$,
and $\lambda(\cdot, t)\in H^{r+1}$. The solution is unique and
depends continuously on the initial data. The conclusions hold
for all $r\ge 0$ if the initial values $\u_0$ and $\u_1$ lie in a 
common plane through the origin. In this case, the solution is planar.
\label{thm:main}
\end{theorem}

Theorem~\ref{thm:main} implies that Eq.~(\ref{eq:EL-x}) with initial values
\begin{equation}
\label{eq:initial-x}
\x(s,0)=\x_0(s)\ ,\quad \partial_t\x(s,0)=\x_1(s)\
\end{equation} 
satisfying  the compatibility condition
\begin{equation} 
\label{eq:compatible-x}
\partial_s\x_0\cdot \partial_s\x_1\equiv 0\ ,
\end{equation}
is well-posed in $H^{7/2}\times H^{3/2}$ (in $H^3\times H^1$, if
the motion is planar).  Periodicity of
$(\x_0,\x_1)$ implies that 
$\x(s+2\pi,t)- \x(s,t)$ and its time derivative are
zero at $t=0$, whereas
\begin{equation}
\label{eq:momentum}
\partial_t^2 \bigl( \x(s+2\pi,t)-\x(s,t)\bigr)
\ =\ \int_0^{2\pi} \partial_t^2 \u(s',t)ds' 
\ =\ 0\ 
\end{equation}
by  Eq.~(\ref{eq:EL-u}), so that the curve $\x(\cdot, t)$ remains
a closed loop.

Theorem~\ref{thm:main} applies also to infinite
(expanding or contracting) helical curves $\x(s,t)$ for
which $\u(\cdot, t)$ is still periodic.
Eq.~(\ref{eq:momentum}) implies that unless $\u_1=\partial_t\u(\cdot, 0)$
averages to zero, $\x(s+2\pi,t)-\x(s,t)$ grows
linearly in time.  Since it cannot grow beyond $2\pi$, the curve 
must disintegrate in a finite time under these circumstances.

Since the total energy ${\cal T}+{\cal V}$ is equivalent to the 
natural norm of $(\x,\partial_t\x)$
in $H^2 \times L^2$, it seems reasonable to 
consider the Eq.~(\ref{eq:EL-x}) in $H^2\times L^2$. 
By conservation of energy, a small-time existence result for solutions of
Eq.~(\ref{eq:EL-x}) in $H^1\times L^2$ would imply that solutions 
exist globally in time.  This would amount to proving
Theorem~\ref{thm:main} with $r=-1$ for initial values $(\u_0,\u_1)$
where $\u_1$ averages to zero.  It is an open question under 
what conditions on the initial values the solutions to Eq.~(\ref{eq:EL-x}) 
exist globally in time.  Blowup, if there is any, must involve a 
transfer of energy to the high-frequency (large $n$) modes of 
the Fourier transform of the solution.

Our theorem should be compared with results of
Caflisch and Maddocks \cite{CM}, who have given
a proof of {\em global} existence for a {\em planar}
dynamical elastica. Their equations of motion include an additional 
rotational inertia term proportional to
$ \frac{1}{2} \int |\partial_t \u|^2 \, ds$
in the kinetic energy ${\cal T}$ (see their Eq.~(2.16)). 
They assume that the initial values
$\u_0$ and $\u_1$ are piecewise $C^2$ and piecewise $C^1$, respectively,
a stronger assumption than we require in the planar case (although it
appears that their smoothness assumptions can be relaxed somewhat).
With this additional kinetic energy term, the equation of motion
(for the tangent  vector angle) can be written as a single
{\em second order} semilinear wave equation with a non-local
nonlinear term. Conservation of energy ultimately
provides the {\em a priori} bounds needed to prove their global
existence result.

\subsection{Related geometrical and physical problems}

The variational problem stated in Eq.~(\ref{lageq}) is
closely related with wave maps.
For a wave map problem, a 
typical choice for the kinetic energy would be
\begin{equation}
{\cal T}_1
\ =\ \frac{1}{2}\int_0^{2\pi}  
             \abs{\partial_t\u(s,t)}^2\, \frac{ds}{2\pi}\ ,
\end{equation}
rather than that of Eq.~(\ref{kineq})
the potential energy is given by Eq.~(\ref{poteq}) with 
$\kappa = |\partial_s\u|$, and $u$ is constrained to
lie in a Riemannian manifold (see \cite{Sh-St}).  
The resulting Euler-Lagrange equations
for $\u$ are just second order in space and time,
and  the Lagrange multipliers can be expressed
in terms of the second fundamental form of
the target manifold applied to first derivatives of $u$.
Both local and global existence results
have been obtained for various wave map 
problems~\cite{Gu,KT,Tat,Tao,KR,NSU}.

We note that choices for the potential energy of a geometric nature 
other than that of Eq.~(\ref{poteq}) are possible. Let 
$\n(\cdot,t)$ and $\b(\cdot,t)$ be the standard normal and the binormal
to the curve $\x(\cdot,t)$, as defined by the Serret-Frenet 
formulas, and let $A=A(s,t)$ be the $3\times 3$-unitary matrix whose columns
are $\u(s,t)$, $\n(s,t)$ and $\b(s,t)$ (see Eq.~(\ref{eq:SF}) below).
An alternate choice for the potential energy related to the
Dirichlet form energies for harmonic maps~\cite{GV,Sh,ER,H,LET} 
\begin{eqnarray}
\label{eq:def-V1}
    {\cal V}_1&=&
    \frac{1}{4}\int_0^{2\pi}tr\left(\partial_s
        A^\dagger(s,t)\,\partial_s A(s,t)\right)\,\frac{ds}{2\pi}\nonumber\\
    &=& \frac{1}{2}\int_0^{2\pi}
         \bigl(\kappa^2(s,t)+\theta^2(s,t)\bigr)\,\frac{ds}{2\pi}
    \end{eqnarray}
where $\kappa$ is the curvature, and $\theta$ is the torsion of the curve. 
Combining this with the kinetic energy in Eq.~(\ref{kineq})
results in  adding  a term of the form
\begin{equation}
   \partial_s^2  \left(\frac{\theta}{\kappa^2} 
          \partial_s\u \times\partial_s^2\u\right)
 + \partial_s^3\left(\frac{\partial_s\theta}{\kappa^2}
          \partial_s\u \times\partial_s^2\u\right) 
\end{equation}
to the right hand side of Eq.~(\ref{eq:EL-u}),
which is $6^{th}$ order in the spatial derivatives of $\u$.  

For a physical closed loop of wire with a small circular
cross section of radius $\rho>0$, a more realistic 
expression for the elastic potential energy is given by
\begin{equation}
\label{eq:def-V2}
{\cal V}_2
\ =\ \frac{\pi}{2} \rho^4 \left\{ (\frac{\lambda}{4}+\frac{\mu}{2})
   \int_0^{2\pi} \! \kappa(s,t)^2 \,\frac{ds}{2\pi}
+ \frac{\mu}{2} \Bigl[ \int_0^{2\pi} \! \theta(s,t)\, 
               \frac{ds}{2\pi}\Bigr]^2\right\}
+ o(\rho^4)\ ,
\end{equation}
where $\lambda$ and $\mu$ are the Lam\' e constants 
for a homogeneous isotropic hyperelastic material
(see~\cite{Ciarlet}).  (For a wire made of material of a fixed density,
the kinetic energy is of order $\rho^2$. The kinetic 
energy due to twisting is of order $\omega^2\rho^4$, where
$\omega$ is the angular velocity.)
The torsion term accounts for the expense of twisting the material 
frame. The equation is derived under the ``quasi-static'' assumption that
the material arranges itself instantaneously
about the central curve such that the contribution 
of the local twist to the elastic energy is
as small as possible. The minimizing configuration of the material
for a given curve is achieved by a local twist which is constant 
along the wire. For a dynamical elastica, that constant will in 
general change over time.  The corresponding potential energy for 
a very thin, narrow ribbon would contain a torsion term proportional 
to the last term in Eq.~(\ref{eq:def-V1}).
Maddocks and Dichmann \cite{MD}, Coleman et al. \cite{CDLT} and
others consider {\em director} theories, originated by Kirchhoff
and Clebsch, in which there are further stress-strain relations 
between the tangent and two independent normal vectors (see \cite{A}).
The above choices for the potential and kinetic energy
will not be pursued here.

Our equation is related to the Localized Induction Equation (LIE)
first discussed by Da Rios \cite{daR}
and rediscovered by Arms and Hama \cite{AH}. This
is a second order equation for the
approximate time evolution of a thin vortex filament $\x(s,t)$ (e.g. a
smoke-ring) moving in a fluid, and is given by
\begin{equation}
     \partial_t\x(s,t)= \partial_s\x\times\partial_s^2\x(s,t),
\label{eq:LIE}
\end{equation}
again, $s$ being arclength. (The equation ignores long-distance
effects of the vortex acting upon itself.)  Differentiating this
equation with respect to $s$,
one obtains an equation for the tangent vector $\u$ to the vortex,
\begin{equation}
     \partial_t\u(s,t)= \u\times\partial_s^2\u(s,t),
\label{eq:LL}
\end{equation}
known as the Landau-Lifshitz equation for the continuum Heisenberg
ferromagnet~\cite{Faddeev}.  We will elaborate on the connections
between our equations and particularly the Landau-Lifshitz equation in
the next subsection.  A convenient reference for the LIE equation can be
found in Newton \cite{Newton}.

The equations of motion for the elastica exhibit a rich 
variety of solutions. Langer and Singer found a
countably infinite number of equilibrium configurations 
which are contained in tori of revolution and represent 
all but a finite number of torus knots, and showed that, 
up to the symmetries of the equations, every non-planar equilibrium 
configuration appears in their list ~\cite{LS1}.  
They also relate these solutions to those of LIE~\cite{LS2}.

At least in modern times, consideration of small
amplitude vibrations of a rotating (in general extensible)
ring seem to go back to Carrier's 1945 paper~\cite{Car}. He moreover
considered the case in which the ring was constrained,
or supported at points around the ring.
Simmonds \cite{S} also considered
small planar vibrational modes for a nearly circular, extensible
ring, in particular flexing modes
in which extension is essentially negligible. His analysis 
provides a systematic treatment of the small amplitude approximations, 
and he shows for example that the vibrational frequencies decrease with 
amplitude.  In a different direction, Coleman and Dill \cite{CD} 
examined the infinite length planar elastica, and showed that the solitary
waves are of the form of a single loop and that
traveling waves are a succession of periodically spaced loops, 
all of which satisfy differential equations similar to Euler's 
equations for the static case.  Note that their equations can 
include rotational inertia, cf. their equations (44a,b). 
They also find wave solutions that are periodic in time.   
Following this work, Coleman and Xu \cite{CX} numerically 
investigated solitary waves for an elastica of large length, 
and showed that the scattering was more than a simple phase shift, thereby
providing compelling evidence that the elastica is not completely
integrable. In \cite{CDLT}, Coleman et al. considered the
elastica moving in ${\bf R}^3$, showing existence of traveling
and solitary waves exhibiting torsion so that the resulting curves
are corkscrew-like. They include some discussion of the (finite)
closed curve case.
In their work on global existence of solutions for the planar elastica, Caflisch and Maddocks \cite{CM} also showed Liapunov stability for
solutions near isolated relative minima of the potential ${\cal V}$.

\subsection{Description of the proof}

A first step towards solving Eq.~(\ref{eq:EL-u})
is to obtain an equation for the Lagrange multiplier~$\lambda$ 
in terms of the solution $\u$.  In the related case of a wave map, 
an explicit expression
for the Lagrange multiplier in terms of first derivatives 
of the solution is obtained by
projecting the equation onto the normal of the target sphere~\cite{Sh-St}.
In our problem, taking the inner product of Eq.~(\ref{eq:EL-u}) with $\u(s,t)$, 
writing $\abs{\partial_s\u(s,t)}=\kappa(s,t)$, and
using that $\u\cdot\partial_s\,\u= \u\cdot\partial_t\,\u=0$
due to the constraint yields an elliptic boundary value problem 
for the tension $\lambda(\cdot,t)$ at time $t$,
\begin{equation}
 \left(-\partial_s^2+\kappa^2\right)\lambda \ =\  \abs{\partial_t\u}^2-
   \u\cdot\partial_s^4\u\ .
\end{equation}
A slightly less obnoxious form of this equation 
(involving lower order derivatives) results if one uses the 
identity 
\begin{equation}\label{eq:bianchi}
    3\abs{\partial_s^2\u}^2 + 4\partial_s\u\cdot\partial_s^3\u
  + \u\cdot\partial_s^4\u \  = \ 0 \  ,
\end{equation}
which follows by differentiating the constraint
$\abs{\u}^2=1$ four times with respect to $s$. With this 
identity the tension equation can be written as 
\begin{equation}
\label{eq:lambda}
\left(- \partial_s^2+\kappa^2\right)
\left(\lambda + 2\kappa^2\right)
\ =\ \abs{\partial_t\u}^2 + 2 \kappa^4 - \abs{\partial_s^2\u}^2\  .
\end{equation}
We do not know how to make sense of Eq.~(\ref{eq:lambda}) without 
requiring at least $\u(\cdot, t)\in H^2$, $\partial_t\u(\cdot, t)\in L^2$.

Rewriting Eq.~(\ref{eq:EL-u}) as a system, we want to solve 
\begin{equation}
\label{eq:system}
\left\{
\begin{array}{rcl}
\partial_t\u &=& \v\\
\partial_t\v &=& -\partial_s^4 \u + \partial_s^2(\lambda\u)\ ,
\end{array}\right.
\end{equation}
where $\lambda$ is determined by Eq.~(\ref{eq:lambda}).
The linear part of Eq.~(\ref{eq:system}) generates  a strongly continuous 
semigroup on $H^{r+2}\times H^{r}$ for any $r\in\R$. Standard 
semilinear theory accommodates a nonlinearity that defines a 
locally Lipschitz continuous map from this space to itself as a perturbation.
Since $\u\in H^{r+2}$, it would suffice to show that
$\lambda\in H^{r+2}$ to apply this technique. However, 
we only have $\kappa\in H^{r+1}$, hence $\kappa^2\in H^{r+1}$
(provided $r>-1/2$ ensuring that $H^{r+1}$ is an algebra)
while $\lambda+2\kappa^2$ is more regular by Lemma~\ref{lem:mu} .
Thus, we can only expect $\lambda\in H^{r+1}$
and so Eq.~(\ref{eq:EL-u}) cannot be solved directly
just by converting it to a  Duhamel integral representation.  
A similar picture emerges if $\lambda$ is inserted into Eq.~(\ref{eq:EL-x}). 

Instead, we proceed as follows.  Since Eq.~(\ref{eq:lambda}) is the projection
of Eq.~(\ref{eq:EL-u}) in the direction of $\u$, 
we combine it with the complementary projection 
\begin{equation}
\begin{array}{c}
-\u\times (\u\times \partial_t^2\u) 
\ =\ \u\times (\u\times \partial_s^4\u)
+ 2(\partial_s \lambda)\partial_s\u 
- \lambda \u\times (\u\times \partial_s^2\u)\ .\end{array}
\label{eq:LL2}
\end{equation}
We have used that $\u\times\u=0$ and $\u\cdot \partial_s\u=0$.
The last two terms on the right hand side are locally 
Lipschitz from $H^{r+2}\times H^r$ to $H^r$.
Now the problem is that the projection of the fourth derivative 
onto the orthogonal complement of $\u$ is a complicated nonlinear 
operation.  
Eqn.~(\ref{eq:LL2}) is related to the Landau-Lifshitz
equation, Eq.~(\ref{eq:LL}). This becomes apparent if we
use Eq.~(\ref{eq:LL}) to formally write a differential equation
for the {\em second} time derivative of $\u$.  The resulting
equation differs from our Eq.~(\ref{eq:LL2}) in the
last two terms: for us the tension is given by the solution to an
elliptic problem, hence these terms are not local,
while for the Landau-Lifshitz case the tension is
given by a local expression of $\u$ and its derivatives. The
similarity between the equations suggests how to proceed.

The Landau-Lifshitz equation 
is equivalent to a nonlinear Schr\"odinger equation via the so-called
{\em Hasimoto transformation} (\cite{Has}, cf.~\cite{CSU}),
which is is a nonlinear, solution-dependent 
change of variables where all partial derivatives of 
a curve are expressed with respect to a local coordinate frame 
transported along the curve.  In Section~\ref{sect:ChVar}, we will
use the Hasimoto transformation
to transform the projection of the second derivative operator
on the left hand side and the fourth derivative operator
on the right hand side of Eq.~(\ref{eq:LL2}) into linear differential
operators plus a perturbation, while leaving the
other terms essentially invariant (see Eq.~(\ref{eq:pq-all})).
Unless the motion described by Eq.~(\ref{eq:EL-u}) is planar, the 
Hasimoto transformation introduces
a {\em monodromy} into the problem. In other words, 
the frame transported along the elastica is not periodic 
but rather quasiperiodic, with a generally {\em time-dependent} rotation
of phase $2\pi \beta(t)$ about the tangent
equal to the integral of the torsion over the curve
(see Eq.~(\ref{eq:def-beta})). We correct
for the monodromy by performing an additional
gauge transformation.  Interestingly, the 
Hasimoto coordinate frame with the monodromy correction
coincides in our problem with the natural material 
frame~\cite{LS2}.  
However, monodromy does not arise 
in the Kirchhoff-Clebsch director theories \cite{MD,CDLT}.
Our Hasimoto transformation allows us to 
rewrite Eq.~(\ref{eq:EL-u}) in the form
\begin{equation}
\frac{d}{dt} Y(t)\ =\ G_{\beta(t)}\,Y(t) \ +\ F_{\beta(t)}(Y(t))\ ,
\label{eq:quasilinear-Y}
\end{equation}
where the linear part $G_{\beta(t)}$ is a differential operator
which generates a strongly continuous evolution operator 
on a Sobolev space $\Y^r$ of periodic functions,
and $F_{\beta(t)}$ is a nonlinear perturbation.
The operator $G_{\beta(t)}$ depends
on the on the monodromy, which is 
in turn determined by an auxiliary equation
\begin{equation}
\label{eq:quasilinear-beta}
\frac{d}{dt}\beta(t)\ =\ B(Y(t))\   ,
\end{equation}
see Eqs.~(\ref{eq:PQ-all}) and (\ref{eq:Y-Duhamel}).

In Section~\ref{sect:estimates}, we provide 
the basic estimates for the linear evolution operator $V_\beta$ generated by
$G_{\beta(t)}$ and for the nonlinearity $F_{\beta(t)}$.
In Lemma~\ref{lem:spectrum},
we prove a positive lower bound on the spectrum
of the Schr\"odinger operator $L_\kappa=-\partial_s^2+\kappa^2$ 
appearing on the left hand side of Eq.~(\ref{eq:lambda}) which may be of 
independent interest.   We conjecture  that the lowest eigenvalue
of $L_\kappa$ is minimal when the elastica is a circle
(see \cite{HL} for a related result on $-\partial_s^2 -\kappa^2$).

In Section~\ref{sect:local}, we set up a contraction mapping argument 
for Eqs.~(\ref{eq:quasilinear-Y})-(\ref{eq:quasilinear-beta}). 
A technical complication is that the dependence
of the linear evolution operator $V_\beta$ generated by $G_{\beta(t)}$
is only strongly continuous, not norm-continuous
in $\Y^r\times\R$.  We overcome this difficulty by setting up
the contraction mapping argument using a weaker norm.

Our local existence proof does not incorporate the more 
modern space-time methods such as Strichartz inequalities.
Versions of these of these inequalities 
adapted to problems with periodic boundary 
conditions~\cite{B} have been used to obtain
global existence of solutions for other related wave equations,
for example the Boussinesq equation on a circle,
which is also fourth order~\cite{FG}.
The methods do not immediately  apply here however.
The (spatial) derivative in the nonlinear term
of the Boussinesq equation is of lower
order than in our wave equation Eq.~(\ref{eq:EL-u}). Thus
the Duhamel integral form of the equations of
motion analogous to our Eqs.~(\ref{eq:system})
make sense for example in $H^{r+1}\times H^r$
for suitable $r$, unlike the situation here which necessitates the Hasimoto
change of variables. But again this change of variables
leads to the monodromy issue and an associated
linear problem with a {\em time-dependent}
generator 
for which we do  not have suitable space-time
estimates (see Eq.~(\ref{eq:quasilinear-Y})). 
Finding such estimates could be an  avenue 
to relaxing our regularity assumptions  for the initial data.

For the dynamical elastica  moving in three dimensions,
the question of global existence remains open. The picture which does
emerge from our approach is that we can always integrate forward for
an open interval of time, i.e. we have existence (and even uniqueness
for planar motion), until 
$\norm{(\u(\cdot, t),\partial_t\u(\cdot, t))}_{H^{2}\times L^{2}}$ 
becomes unbounded.  
But at this moment, the tension $\lambda(s)$ given by
Eqs.~(\ref{eq:lambda}) becomes infinite at some point $s$ 
(or at least the individual terms on the right hand side of 
Eq.~(\ref{eq:lambda}) are not integrable
so that it is by no means clear
that at this moment $\lambda$ exists even as a distribution).
Seemingly the elastica would break apart. 
It would be of interest to know whether indeed infinite tension can
develop in a finite time. 

We conclude with a couple of remarks about the infinite length
elastica. A lower bound on the spectrum
of $L_\kappa$ is by no means apparent in this case,
for example $L_{\kappa}$ acting in $L^{2}(\mathbf{R})$ is 
typically {\em not} invertible, i.e. there is an infrared 
divergence and the estimates in Lemma~\ref{lem:spectrum} would fail. 
Physically, this divergence corresponds to
elastica configurations with infinite tension. 
A local existence proof for the infinite elastica would thus 
presumably involve more  subtly defined function spaces for 
which the tension is finite.
The second question concerns the role of the monodromy introduced
by the Hasimoto transformation. Even for finite-energy solutions, it
is not obvious how to remove the total twist by a gauge transformation.  

\section{A Change of Variables}
\label{sect:ChVar}
\setcounter{equation}{0}

\subsection{Hasimoto transformation}

Let $\u(s,t)$ be a smooth solution of Eq.~(\ref{eq:EL-u}). 
We will express the partial derivatives of $\u$ in terms
of a positively oriented orthonormal frame which
consists of $\u$ and two other unit vectors 
which we combine to a single complex vector
$\tilde\v$. 
The vector $\tilde \v(s,t)$ is chosen so that for any fixed time $t$,
\begin{equation}
\partial_s \tilde \v\ =\ -(\tilde \v\cdot \partial_s\u)\, \u\ ,
\label{eq:parallel}
\end{equation}
that is, the real and imaginary parts of
$\tilde \v(\cdot, t)$ are moved along the
curve $\u(\cdot, t)$ by {\em parallel transport} on $\S^2$.  Then
\begin{equation}
\left\{
\begin{array}{lcl}
\partial_s \u&=&\Re{\bar q \tilde \v}\\
\partial_s \tilde \v&=&-q\u
\end{array}\right.
\qquad\qquad
\left\{
\begin{array}{lcl}
\partial_t \u&=&\Re{\bar p \tilde \v}\\
\partial_t \tilde \v&=&-p\u + i\tilde \alpha \tilde \v
\end{array}
\right.
\label{eq:uv0-ODE}
\end{equation}
where $p(s,t)$ and $q(s,t)$ are complex-valued, and
$\tilde \alpha(s,t)$ is real-valued.
Since $\partial_s\partial_t \u = \partial_t\partial_s \u$
and $\partial_s\partial_t \tilde \v = \partial_t\partial_s\tilde \v$, 
we must have
\begin{equation}
(\partial_t-i\tilde \alpha) q \ =\ \partial_s p, \quad
\partial_s\tilde \alpha = \Im{\bar q p }.
\label{eq:consistent}
\end{equation}
In these coordinates, we compute for
the projection of Eq.~(\ref{eq:EL-u}) in the direction of~$\u$,
\begin{equation}
-\abs{p}^2\ =\ \partial_s^2 \left(2\abs{q}^2 +\lambda\right)
- \abs{q}^2\left(\abs{q}^2 +\lambda \right) - \abs{\partial_s q}^2\ ,
\label{eq:lambda-1}
\end{equation}
see Eq.~(\ref{eq:lambda}). For the complementary projection, we find 
\begin{eqnarray}
\label{eq:LL-pq}
(\partial_t -i\tilde \alpha)p 
&=& \partial_s^2\left(-\partial_sq\right) 
    +2 \partial_s\left((\abs{q}^2 +\lambda) q\right) -
\left(\abs{q}^2+\lambda\right)\partial_sq +\Re{\bar q \partial_s q}q\ ,\quad
\end{eqnarray}
see Eq.~(\ref{eq:LL2}).  If $\u$ is a smooth solution of
Eq.~(\ref{eq:EL-u}), and $\tilde \v$, $\tilde \alpha$, $p$, $q$,
are defined by Eqs.~(\ref{eq:uv0-ODE}), and if we set
\begin{equation}
\mu\ =\ \lambda + 2\abs{q}^2\ ,
\label{eq:lambda-mu}
\end{equation}
where $\kappa=\abs{q}$ is the curvature, then we arrive at the system
\begin{equation}
\left\{
\begin{array}{rcl}
(\partial_t -i\tilde \alpha)p 
  &=& - \partial_s^3 q -4 \Re{\bar q\, \partial_sq}q
      - i \Im {\bar q\, \partial_s q}q \\ 
 && \qquad + 2(\partial_s\mu)q +\mu\partial_s q\\
(\partial_t -i\tilde \alpha) q 
  &=& \partial_s p \\
\partial_s \tilde \alpha 
  &=& \Im{\bar q p} \\
(-\partial_s^2+ \abs{q}^2)\mu
  &=& \abs{p}^2 + \abs{q}^4- \abs{\partial_s q}^2 \ .
\end{array}\right.
\label{eq:pq-all}
\end{equation}
Here, the first equation follows from Eq.~(\ref{eq:LL-pq}), 
the second and third are the consistency relations in 
Eq.~(\ref{eq:consistent}), and the last follows from 
Eqs.~(\ref{eq:lambda-1}) and~(\ref{eq:lambda-mu}).

\begin{remark} {\em 
There are many other ways to complete $\u$ to an orthonormal 
frame, which are all related by {\em gauge transformations}
\begin{equation}
\begin{array}{c}
\tilde \v \mapsto  e^{i\gamma(s,t)} \tilde \v\ ,\quad
p \mapsto e^{i\gamma(s,t)} p\ ,\quad
q \mapsto e^{i\gamma(s,t)} q\ ,  \\
\tilde \alpha \mapsto \partial_t\gamma(s,t)+ \tilde \alpha\ ,\quad
\mu \mapsto \mu 
\end{array}
\label{eq:gauge}
\end{equation}
with some real-valued function $\gamma(s,t)$.
The choice in Eq.~(\ref{eq:parallel}) has 
the special property that no second derivatives
appear on the right hand sides of Eq.~(\ref{eq:pq-all}).
As discussed in connection with Eq.~(\ref{eq:LL2}), this is a key step 
towards solving Eqs.~(\ref{eq:EL-u})-(\ref{eq:initial-u}),
because expressions containing first derivatives of $q$
can, but expressions containing second derivatives
of $q$ {\em cannot} be treated as perturbations of the 
third derivative operator in the first line of Eq.~(\ref{eq:pq-all}).
}\end{remark}

\subsection{The monodromy correction}

Even when $\x(\cdot,t)$ is $2\pi$-periodic, the 
frame $(\u,\tilde \v)$ defined by Eq.~(\ref{eq:parallel})
is in general quasiperiodic,
\begin{equation}
\tilde \v(s\!+\!2\pi, t) \ =\ e^{2\pi i\beta(t)} \tilde \v(s,t)\ ,
\end{equation}
where $\beta$ is a real-valued function  of time.
The value of $\beta(t)$ is determined modulo an integer 
by the {\em monodromy} of the parallel transported frame, which is
related to the curve torsion, see Eq.~(\ref{eq:def-beta}).
It follows that $p$, $q$, and $\alpha$  satisfy 
time-dependent quasiperiodic boundary conditions.

Periodic boundary conditions are recovered by the gauge transformation
\begin{equation}
\begin{array}{rcl}
\v(s,t)&=& e^{i(\tilde\beta(t)-s\beta(t))} \tilde \v(s,t)\\
P(s,t)&=& e^{i(\tilde\beta(t)-s\beta(t))} p(s,t)\\
Q(s,t)&=& e^{i(\tilde\beta(t)-s\beta(t))} q(s,t)\\
\alpha(s,t) &=& \tilde \alpha(s,t) +
\partial_t\bigl(\tilde\beta(t)-s\beta(t)\bigr) \ ,
\end{array}
\label{eq:pq-PQ}
\end{equation}
where $\beta(t)$ and $\tilde \beta(t)$ are chosen so that
$\v$ is periodic, and $\alpha$ averages to zero
over a period of~$2\pi$.  We arrive at the system of equations
\begin{equation}
\left\{
\begin{array}{rcl}
(\partial_t -i\alpha) P 
 &=& -(\partial_s+i\beta)^3 Q 
       -4 \Re{\bar Q\,\partial_s Q}Q- i\Im{\bar Q\,\partial_sQ}Q\\
 && \qquad - i\beta \abs{Q}^2Q
+ 2 (\partial_s\mu)Q +\mu \,(\partial_s\!+\!i\beta)Q\\
(\partial_t -i\alpha) Q 
 &=& (\partial_s+i\beta) P \\[0.1cm]
\partial_s \alpha +\partial_t\beta
 &=& \Im {\bar Q P} \\[0.1cm]
\left(-\partial_s^2+\abs{Q}^2\right)\mu 
 &=& \abs{P}^2 +\abs{Q}^4 - \abs{(\partial_s\!+\!i\beta)Q}^2\ .
\end{array}\right.
\label{eq:PQ-all}
\end{equation}
Here, $P(s,t), Q(s,t)\in\C$, and $\alpha(s,t), \beta(t), \mu(s,t)\in \R$
satisfy $2\pi$-periodic boundary conditions, with the constraint
\begin{equation}
\int_0^{2\pi}\alpha(s,t)\,\frac{ds}{2\pi} \ =\ 0\ .
\label{eq:mean0}
\end{equation}
Initial conditions are given by
\begin{equation}
P(s,0)\ =\ P_0(s)\ ,\quad Q(s,0)\ =\ Q_0(s)\ ,\quad \beta(0)=\beta_0\ ,
\label{eq:initial-PQbeta}
\end{equation}
where $P_0$ and $Q_0$ are periodic complex-valued functions, and
$\beta_0$ is the monodromy at time $t=0$. 
We emphasize here that  it is really this system,
Eqs.~(\ref{eq:PQ-all})-(\ref{eq:initial-PQbeta}) that we analyze in 
this paper.  We will see in Section~\ref{sect:estimates} 
that Eq.~(\ref{eq:PQ-all}) can be written in the form of 
Eqs.~(\ref{eq:quasilinear-Y})-(\ref{eq:quasilinear-beta})
with $Y=(P,Q)$.

\begin{remark} {\em Eqs.~(\ref{eq:PQ-all})-(\ref{eq:mean0})
and the periodic boundary conditions are invariant under the gauge 
transformations
\begin{equation}
\begin{array}{c}
P\mapsto e^{i(s_0+ks)}P\ ,\quad
Q\mapsto e^{i(s_0+ks)}Q\ ,\\
\alpha\mapsto\alpha\ ,\quad
\beta \mapsto \beta - k\ ,\quad 
\mu\mapsto\mu\ ,
\end{array}
\label{eq:PQ-gauge}
\end{equation}
when $k$ is an integer and $s_0\in\R\,$.
}\end{remark}

We turn to the relation of the initial-value problem in
Eqs.~(\ref{eq:PQ-all})-(\ref{eq:initial-PQbeta})
with the initial-value problem in Eqs.~(\ref{eq:EL-u})-(\ref{eq:initial-u})
posed in the introduction.  We  say the initial conditions 
$(P_0,Q_0,\beta_0)$ for Eq.~(\ref{eq:PQ-all}) are {\em compatible with} 
Eq.~(\ref{eq:EL-u}), if the linear system
\begin{equation}
\left\{
\begin{array}{lcl}
\partial_s \u&=&\Re{\bar Q_0 \v}\\
\partial_s \v&=&-Q_0\u-i\beta_0 \v\ ,
\end{array}\right.
\label{eq:uv-ODE-1}
\end{equation}
has a $2\pi$-periodic solution forming an orthonormal frame.

\begin{lemma} \label{lem:equiv}
For $r\ge 0$, the following statements are equivalent:
\begin{enumerate}

\item For each pair of initial values $(\u_0,\u_1)\in H^{r+2}\times H^r$ 
satisfying the condition in Eq.~(\ref{eq:compatible-u}), the 
initial value problem in Eqs.~(\ref{eq:EL-u})-(\ref{eq:initial-u})
has a solution $(\u,\partial_t\u)\in H^{r+2}\times H^r$,
defined on some short time interval, which is strongly
continuous in $t$ and assumes the initial values at $t=0$.
The solution is unique and depends continuously on the initial values.
If the initial values are smooth then the solution is smooth
in both variables.

\item For each triple of initial values 
$(P_0,Q_0,\beta_0)\in H^{r}\times H^{r+1}\times\R$
which is compatible with Eq.~(\ref{eq:EL-u}),
the initial value problem in
Eqs.~(\ref{eq:PQ-all})-(\ref{eq:initial-PQbeta}) has a solution 
$(P,Q,\beta)\in H^r\times H^{r+1}\times \R$,
defined on a short time interval, which is 
strongly continuous in $t$ and assumes the initial values at $t=0$.
The solution is unique and depends continuously on the initial values.
If the initial values are smooth then the solution is smooth
in both variables.

\end{enumerate}
\end{lemma}

\begin{remark} {\em  The regularity assumption  $r\ge 0$ 
is needed only to ensure that the right hand side of the equation for
$\mu$ in the fourth line of Eq.~(\ref{eq:PQ-all}) makes sense
as an $L^1$-function.  In the proof of the lemma,
we will construct a transformation $(\u,\partial_t \u) \mapsto (P,Q,\beta)$ 
between $H^{r+2}\times H^r $ and $ H^r\times H^{r+1}\times \R$
which is continuous and has a continuous inverse
for all $r\ge -1$. 
}
\end{remark}

\medskip\noindent{\sc Proof of Lemma~\ref{lem:equiv}: \ }\nobreak
Given initial values $(\u_0,\u_1)$ for Eq.~(\ref{eq:EL-u})
satisfying the condition in Eq.~(\ref{eq:compatible-u}), we determine
initial values $(P_0,Q_0,\beta_0)$ for Eq.~(\ref{eq:PQ-all})
by choosing a complex vector $\tilde \v_0(0)$ whose 
real and imaginary parts complement
$\u_0(0)$ to a positively oriented orthonormal basis of $\R^3$,
and then solving Eq.~(\ref{eq:parallel}) with $\u(s,0)=\u_0(s)$
to obtain $\tilde \v_0(s)$.
The initial monodromy $\beta_0$ is determined up to an additive 
integer by
\begin{equation}
\tilde \v_0(2\pi)\ =\ e^{2\pi i \beta_0}\tilde \v_0(0)\ ,
\label{eq:init-beta}
\end{equation} 
and $P_0$ and $Q_0$ are given by
\begin{equation}
P_0(s)\ =\ e^{-is\beta_0}\tilde \v_0(s)\cdot \u_1(s)\ ,\quad
Q_0(s)\ =\ e^{-is\beta_0} \tilde \v_0(s)\cdot \partial_s \u_0(s)\ .
\label{eq:init-pq}
\end{equation}
This defines a continuous map from
$H^{r+2}\times H^{r}$ to $H^r\times H^{r+1}\times\R$
for any $r\ge 0$. 
We have seen above that this transformation maps
smooth solutions $u$ of Eq.~(\ref{eq:EL-u})-(\ref{eq:BC-u}) 
to smooth solutions $(P,Q,\beta)$ of Eq.~(\ref{eq:PQ-all})-(\ref{eq:mean0}).
Since smooth functions are dense in $H^{r+2}\times H^r$,
the transformation can be extended continuously to
all of $H^{r+2}\times H^r$.

Conversely, given $(P_0,Q_0,\beta_0)$, let $(\u_0,\v_0)$ be 
a solution of the linear differential equation 
Eq.~(\ref{eq:uv-ODE-1}) which defines a 
periodic orthonormal frame, and set
\begin{equation}
\u_1\ =\ \Re{\bar P_0 \v_0}\ .
\end{equation}
By construction, $(\u_0,\u_1)$ satisfy Eq.~(\ref{eq:compatible-u}).
This defines a continuous transformation
from $H^r\times H^{r+1}\times \R$ to $H^{r+2}\times H^r$.
Given a smooth solution of Eqs.~(\ref{eq:PQ-all})-(\ref{eq:mean0}),
we can define a frame $(\u,\v)$ by solving
\begin{equation}
\left\{
\begin{array}{lcl}
\partial_s \u&=&\Re{\bar Q \v}\\
\partial_s \v&=&-Q\u-i\beta \v\ ,
\end{array}\right.
\qquad
\left\{
\begin{array}{lcl}
\partial_t \u&=&\Re{\bar P \v}\\
\partial_t \v&=&-P\u-i\alpha \v\ 
\end{array}\right.
\label{eq:uv-ODE}
\end{equation}
with initial conditions $\u(s,0)=\u_0(s)$, $\v(s,0)=\v_0(s)$.
Note that the third equation in Eq.~(\ref{eq:PQ-all})
ensures that the two systems in Eq.~(\ref{eq:uv-ODE})
can be solved simultaneously. Since $\u_0$ and $\v_0$ are periodic,
$\u(\cdot,t)$ and $\v(\cdot,t)$ are periodic for $t>0$ because 
the pair of equations on the right hand side of Eq.~(\ref{eq:uv-ODE}) 
preserves periodicity by the periodicity of $P$ and $\alpha$.
The function $\u$ obtained in this way is smooth and solves
Eqs.~(\ref{eq:EL-u})-(\ref{eq:initial-u}).
As above, the transformation can be extended
continuously from the subset of
smooth functions to all of $H^r\times H^{r+1}\times\R\,$.

\hfill $\Box$

\subsection{The standard normal frame} 

It is instructive to express $\v$, $\alpha$, $P$, and $Q$ in
terms of standard normal coordinates, which are defined for
smooth curves at any point where the curvature
does not vanish.  Assume $\u(s,t)$ describes the unit tangent 
vector of such a curve $\x(s,t)$, let $\n(s,t)$ be the unit 
normal to the curve in the direction of $\partial_s\u(s,t)$, 
and set $\b= \u \times \n$. The vector $\n$ is
called the {\em principal normal}, and $\b$ the {\em binormal}
of the curve at $s$, at time $t$. The standard normal frame $(\u,\n,\b)$
is characterized by the Serret-Frenet differential equations
\begin{equation}
\begin{array}{lcl}
\partial_s \u& =& \kappa \n \\
\partial_s (\n+i\b)&=& -\kappa \u -i\theta (\n + i\b)
\end{array}
\label{eq:SF}
\end{equation}
where $\kappa(s,t)$ is the {\em curvature}, and 
$\theta(s,t)$ is the {\em torsion} of the curve.
By definition, $\kappa$ is nonnegative, and $\theta$ is real-valued.
The curvature can be expressed in the various frames as
\begin{equation}
\kappa \ =\ \abs{\partial_s \u}
       \ =\ \n\cdot \partial_s \u\ 
       \ =\ \abs{q}
       \ =\ \abs{Q}\ ,
\label{eq:def-kappa}
\end{equation}
and the torsion as
\begin{equation}
\theta \ =\ \frac{\u\cdot(\partial_s \u \times \partial_s^2\u)}
             {\abs{\partial_s\u}^2}
       \ =\ \b\cdot\partial_s\n\ 
       \ =\ \Im{q^{-1}\,\partial_s q}
       \ =\ \Im{Q^{-1}\,\partial_s Q} + \beta\ .
\label{eq:def-theta}
\end{equation}
If we set
\begin{equation}
\gamma(s,t) \ =\  \int_0^s \theta(s', t)\, ds'
\label{eq:def-gamma}
\end{equation}
and define the functions $\tilde \v(s,t)$, $\tilde \alpha(s,t)$, 
$p(s,t)$, and $q(s,t)$ by
\begin{equation}
\begin{array}{lcl}
\tilde \v &=& e^{i\gamma} (\n+i\b)\\
\tilde \alpha &=& \partial_t\gamma + \n \cdot \partial_t \b
\end{array}
\qquad
\begin{array}{lcl}
p &=& e^{i\gamma} \partial_t \u\cdot (\n+i\b)\\
q &=& e^{i\gamma} \kappa\ ,
\end{array}
\label{eq:unb-vpqalpha}
\end{equation}
then the differential equations in Eq.~(\ref{eq:uv0-ODE}) and
the consistency conditions in Eq.~(\ref{eq:consistent})
are satisfied.  The monodromy $\beta$ can be defined by
\begin{equation}
\beta(t)\ =\ \int_{0}^{2\pi}\theta(s, t)\, \frac{ds}{2\pi}\ ;
\label{eq:def-beta}
\end{equation}
in particular, the monodromy is an intrinsic
quantity associated with the curve. 

\medskip
For planar curves, the Hasimoto transformation and the 
analysis of the transformed initial value problem 
simplify considerably. For a curve lying in the
$x_1$-$x_2$-coordinate plane, we choose $\Re{\v}$ to be the 
unit vector obtained by rotating $\u$ counterclockwise through an angle of
$\pi/2$, and $\Im{\v}$ to be the unit vector in the $x_3$-direction,
and set $\tilde \v=\v$. The partial derivatives of $\u$ and $\v$ 
satisfy Eq.~(\ref{eq:uv0-ODE}) with $\tilde \alpha=0$ and $p,q$ 
real-valued. The equations on the left hand side
of Eq.~(\ref{eq:uv0-ODE}) agree with the planar Serret-Frenet equations,
and $q(s,t)$ is the {\em signed curvature} of the curve 
$\x(s, t)$. In this case, no monodromy correction is required
and $(P,Q)=(p,q)$ satisfy Eqs.~(\ref{eq:PQ-all})  with $\alpha=\beta=0$.
The resulting system can be written as a semilinear equation
\begin{equation}
\frac{d}{dt} Y(t)\ =\ G\, Y(t) + F(Y(t))
\end{equation}
in a suitable Sobolev space of periodic functions, which can be solved by
a standard fixed point argument (see Corollary~\ref{cor:planar}).

\section{Estimates}
\label{sect:estimates}
\setcounter{equation}{0}

In the previous section, we have changed
variables in the equations for or the dynamical 
Euler's elastica, and transformed 
Eqs.~(\ref{eq:EL-u})-(\ref{eq:initial-u})
into the equivalent initial value problem given by
Eqs.~(\ref{eq:PQ-all})-(\ref{eq:initial-PQbeta}).
To see that the resulting equations
have the form of Eq.~(\ref{eq:quasilinear-Y})-(\ref{eq:quasilinear-beta}), 
set $Y=(P,Q)$, and define the linear operator $G_{\beta(t)}$ by 
\begin{equation}
G_{\beta(t)} \left(\!\!\begin{array}{c} P\\Q\end{array}\!\!\right)
\ =\ \left(\begin{array}{cc} 0 & -(\partial_s+i\beta(t))^3\\
       (\partial_s+i\beta(t)) & 0 \end{array}\right)
       \left(\!\!\begin{array}{c} P\\Q\end{array}\!\!\right)\ ,
\label{eq:def-G}
\end{equation}
where $\beta$ is is a real-valued function of time.
Suppressing the time-dependence in the notation,
we write the nonlinearity in Eq.~(\ref{eq:quasilinear-Y})
as a sum of three terms:
\begin{eqnarray}
\nonumber
F_\beta(Y)
 &=& i\alpha\left(\!\!\begin{array}{c} P\\Q \end{array}\!\!\right)
+\left(\!\!\begin{array}{c} 2(\partial_s\mu)Q + 
              \mu(\partial_s+i\beta)Q \\ 0\end{array} \!\!\right)\nonumber\\
&& \quad    \ +\ \left(\!\!\begin{array}{c}-4\Re{\bar Q\partial_sQ}Q-
   i\Im{\bar Q\partial_sQ}Q -i\beta \abs{Q}^2Q \\ 0\!\!\end{array}\right)
               \nonumber\\
&=:& F^{(1)}(Y) + F^{(2)}_\beta(Y) + F^{(3)}_\beta(Y)\ , \label{eq:def-F}
\end{eqnarray}
where $\alpha$ is determined by
\begin{equation}
\partial_s \alpha \ =\ \Im{P\bar Q} - 
        \int_0^{2\pi} \Im{P(s)\bar Q(s)}\,\frac{ds}{2\pi}\ ,
\quad \int_0^{2\pi} \alpha(s)\, \frac{ds}{2\pi}=0\ ,
\end{equation}
and $\mu$ solves the elliptic boundary value problem
\begin{equation}
(\abs{Q}^2-\partial_s^2)\mu\ =\ \abs{P}^2 + \abs{Q}^4 
                 - \abs{(\partial_s+i\beta)Q}^2\ 
\label{eq:def-mu}
\end{equation}
with periodic boundary conditions.  The nonlinearity in 
Eq.~(\ref{eq:quasilinear-beta}) is given by
\begin{equation}
B(Y)\ =\ \int_0^{2\pi} \Im{\bar Q(s) P(s)}\, \frac{ds}{2\pi}\ .
               \label{eq:def-B}
\end{equation}

\subsection{The linear part}

We begin the analysis of the system in Eqs.~(\ref{eq:PQ-all})-(\ref{eq:mean0})
by solving the linear equation
\begin{equation}
\frac{d}{dt} \left(\!\!\begin{array}{c} P \\Q \end{array} \!\!\right)
\ =\ G_{\beta(t)} \left( \!\!\begin{array}{c} P\\Q \end{array}\!\!\right)
\label{eq:PQ-linear}
\end{equation}
with a given Lipschitz continuous function $\beta $ defined on $[0,T]$. 
The fundamental solution of Eq.~(\ref{eq:PQ-linear})
will be denoted by $V_\beta(t,t_0)$.  

It is natural to consider Eq.~(\ref{eq:PQ-linear}) in the Fourier
series representation
\begin{eqnarray}\label{eq:Vbeta-1}
\frac{d}{ d t}
       \left(\!\!\begin{array}c \hat P(n,t)\\ \hat Q(n,t) 
                             \end{array}\!\!\right)
   &=& \left(\!\!\begin{array}{cc} 0&i(n+\beta(t))^3\\ i(n+\beta(t))&0 
                             \end{array}\!\!\right)
       \left(\!\!\begin{array}c \hat P (n,t)\\ \hat Q (n,t) 
                            \end{array}\!\!\right)\ ,
\end{eqnarray}
where it decouples into a sequence of linear ordinary differential 
equations on $\C^2$.
Let $\hat G_\beta(n,t)$ be the $2\times 2$ matrix appearing 
on the right hand side of Eq.~(\ref{eq:Vbeta-1}).
The fundamental solution $\hat V_\beta(n,t,t_0)$
of Eq.~(\ref{eq:Vbeta-1}) 
is given by the time-ordered exponential of $\hat G_\beta(n,t)$.

Denote by $H^r$ the Sobolev space of $2\pi$-periodic complex-valued
functions (or distributions, when $r<0$) having $r$
fractional derivatives, with norm
\begin{equation}
\norm{f}_{H^r}^2\ =\ \sum_{n=-\infty}^\infty (1+n^2)^r \abs{\hat f(n)}^2\ .
\end{equation}
Let $\Y^r$ be the space of $2\pi$-periodic functions $Y=(P,Q)$
in $H^r\times H^{r+1}$, with norm
\begin{equation}
\norm{Y}_{\Y^r}^2\ =\ \|P\|^2_{H^r} + \|Q\|^2_{H^{r+1}}\ .
\label{eq:def-norm-Y}
\end{equation}
For a vector $(a,b)$ in ${\bf C}^2$, we define its $n$-norm by
\begin{equation}
\left|\left(\!\! \begin{array}{c} a\\ b
 \end{array}\!\!\right)\right|_n^2
\ =\ |a|^{2}+ w(n)^2 |b|^{2}\ ,
\end{equation}
where $w(n)=\sqrt{1+n^2}$.
With this notation, Eq.~(\ref{eq:def-norm-Y}) becomes
\begin{equation}
\norm{Y}_{\Y^r}^2
\ =\ \sum_{n=-\infty}^\infty w(n)^{2r} \|\hat Y(n)\|_n^2 \ ,
\end{equation}
where $\hat Y(n)=(\hat P(n),\hat Q(n))\in\C^2$ is the Fourier transform
of $Y$.  

\begin{lemma}\label{lem:Vbetabound} Suppose that $\beta$ is 
a real-valued Lipschitz-continuous function
on $[0,T]$ with Lipschitz constant $\eta$, and 
that $\abs{\beta(0)}\le 1$.  There exists an increasing continuous
function $C$ of two variables with $\sup_\eta C(\eta,0)<\infty$ 
such that for $0\le t_0\le t\le T$ and any value of $r\in \R$,
\begin{equation}
\label{claim:Vbetabound}
   \|V_{\beta}(t,t_{0})\|_{\Y^r} \leq C(\eta,T)\ .
\end{equation}
Moreover, $V_\beta(t,t_0)$ is strongly 
continuous in $\Y^r$ with respect to $t$ and $t_0$.
\end{lemma}

\proof
We first show that
\begin{equation}
   \|\hat{V}_{\beta}(n,t,t_{0})\|_n \leq C(\eta,T)\ 
\label{claim:Vbetabound-n}
\end{equation}
uniformly in $n$, with $C$ as in the statement of the lemma.
The idea is to rewrite Eq.~(\ref{eq:Vbeta-1}) as 
\begin{eqnarray}\label{eq:Ubeta-1}
\frac{d}{d t}
       \left(\!\!\begin{array}{c} \hat P(n,t)\\ (n+\beta)\hat Q(n,t) 
                             \end{array}\!\!\right) &&\\
\nonumber && \hskip -2cm 
       =\   i(n+\beta)^2 \left(\!\!\begin{array}{cc} 0&1\\ 1 &0 
                             \end{array}\!\!\right)
       \left(\!\!\begin{array}{c} \hat P (n,t)\\ (n+\beta)\hat Q (n,t) 
                            \end{array}\!\!\right) 
               + \frac{d\beta}{ dt}
          \left(\!\!\begin{array}{c}0\\\hat Q(n,t)\end{array}\!\!\right) \ .
\end{eqnarray}
Let $\hat{U}_\beta(n,t,t_0)$ be the unitary $2\times2$ matrix defined by
\begin{equation}
     \hat{U}_n(\beta,t,t_0)
   \ =\  \exp{\left\{i\int_{t_{0}}^t(n+\beta(t'))^{2}\, dt'\,
        \left(\!\!\begin{array}{cc}
            0&1\\1&0 \end{array}\!\!\right)\right\}} \ .
\label{eq:Ubeta-2}
\end{equation}
By the Duhamel integral formula,
Eq.~(\ref{eq:Ubeta-1}) is equivalent to
\begin{eqnarray}\label{eq:Vbeta-2}
\left(\!\!\begin{array}c
          \hat{P}(n,t)\\
          (n+\beta)\hat{Q}(n,t)
          \end{array}\!\!\right)
    &=& \int_{t_{0}}^t 
                 \hat{U}_{\beta}(n,t,t')
                 \frac{d\beta/dt}{(n+\beta)} 
                 \left(\!\!\begin{array}c 0\\ (n+\beta)\hat{Q}(n,t')
          \end{array}\!\!\right) dt' \nonumber\\
       && \ + \ \hat{U}_{\beta}(n,t,t_{0})\left(\!\!\begin{array}c
          \hat{P}(n,t_{0})\\
         (n+\beta) \hat{Q}(n,t_{0})
          \end{array}\!\!\right)\ 
\end{eqnarray}
provided that $|n|$ is sufficiently large so that $n+\beta(t)$
does not vanish anywhere on $[0,T]$. Since $U_\beta(n, t, t_0)$ is
unitary, and $\abs{\beta(t)}\le 1+\eta T$ for $0\le t\le T$ by assumption,
we can estimate for $|n|\ge 2(1+\eta T)$, 
\begin{eqnarray}
\nonumber 
\left|\left(\!\!\begin{array}c
          \hat{P}(n,t)\\
          (n+\beta)\hat{Q}(n,t)
          \end{array}\!\!\right)\right|
&\le& \frac{2\eta}{|n|} \int_{t_{0}}^t
         \left|\left(\!\!\begin{array}c 0\\ (n+\beta)\hat{Q}(n,t')
                \end{array}\!\!\right)\right| dt' \\
&& \quad + \ \left|\left(\!\!\begin{array}c \hat{P}(n,t_{0})\\
         (n+\beta) \hat{Q}(n,t_{0})
          \end{array}\!\!\right)\right|\ .
\label{eq:Vbeta-3}
\end{eqnarray} 
Applying Gronwall's inequality, and using the fact that the left hand side
is equivalent to the $n$-norm, we arrive at
\begin{equation}
\label{eq:Gronwall-1}
\norm{\hat V_\beta(n,t,t_0)}_n\ \le\  6e^{2\eta T/\abs{n}}\ ,
                       \quad (\abs{n}\ge 1+\eta T)\ .
\end{equation}
For any value of $n$, we can bound the $n$-norm of 
$\hat G_\beta(n,t)$ of Eq.~(\ref{eq:Vbeta-1}) on $\C^2$ by 
\begin{eqnarray}
\norm{\hat G_\beta(n,t)}_n
&=& \sup_{a^2 + w(n)^2b^2=1}
   \left\vert \left(\!\!\begin{array}{cc} 0&i(n+\beta(t))^3\\ i(n+\beta(t))&0 
                     \end{array}\!\!\right)
              \left(\!\!\begin{array}{c}  a\\b
                  \end{array}\!\!\right) \right\vert_n\nonumber\\ 
&\le\ &
\max\left\{ (|n| +1 + \eta T)w(n),\ 
 \frac{(|n|+1 +\eta T)^3}{w(n)}\right\}\label{eq:Gbetabound}\ .
\end{eqnarray}
Estimating the right hand side and applying Gronwall's
inequality gives
\begin{equation}
\label{eq:Gronwall-2}
\norm{\hat V_\beta(n,t,t_0)}_n\ \le\ e^{2(|n| + 1+\eta T)^2(1+\eta T)T}\ .
\end{equation}
Combining Eqs.~(\ref{eq:Gronwall-1}) and (\ref{eq:Gronwall-2})
implies Eq.~({claim:Vbetabound-n}) with the value of the constant given by 
$C(\eta,T)=6 e^{8(1+\eta T)^3 T}$. 
This proves the claim in Eq.~(\ref{claim:Vbetabound}).

Clearly, each Fourier coefficient
$\hat V_\beta(n, t,t_0)\hat Y(n)$ depends continuously on $t$ and 
$t_0$. Since Eq.~(\ref{claim:Vbetabound-n}) implies  a uniform
tail estimate on $\|\hat V_\beta(n,t,t_0)\hat Y(n)\|$, it follows
that $V_\beta(t,t_0)$ is strongly continuous in both
time variables.

\hfill$\Box$

\bigskip\noindent
We also need to bound the dependence of $V_\beta$ on $\beta$.

\begin{lemma} \label{lem:Vbeta-Delta}
Assume that for some $T>0$,
the functions $\beta_1$ and $\beta_2$ are Lipschitz continuous
on $[0,T]$ with Lipschitz constant $\eta$, and that 
$\abs{\beta_1(0)}, \abs{\beta_2(0)}\le 1$.
There exists an increasing continuous function $C$ of two variables 
such that
\begin{equation} 
\label{claim:Vbeta-Delta}
\norm{V_{\beta_2}(t,t_0)-V_{\beta_{1}}(t,t_0)}_{\Y^r \to \Y^{r-1}}
    \ \le\ C(\eta, T)
   \int_{t_0}^t|\beta_2(t')-\beta_{1}(t')|\,dt'\ .
\end{equation}
Moreover, $V_\beta$ is strongly continuous in $\beta$ with 
respect to the $\Y^r$-topology in the sense that for every $Y\in\Y^r$
and any given $\eps>0$ there exists 
$\delta=\delta(\eps,\eta,T,Y)$ such that
\begin{equation}
\sup_{0\le t\le T}\abs{\beta_2(t)-\beta_1(t)}\ \le\ \delta
\end{equation}
implies
\begin{equation}
\sup_{0\le t_0\le t\le T} 
    \norm{(V_{\beta_2}(t,t_0)-V_{\beta_1}(t,t_0))Y}_{\Y^r}\ 
        \ \le\ \eps\ .
\end{equation}
\end{lemma}

\proof We first show that
for $0\leq t_{0}\leq t\leq T$,
\begin{equation}
\label{claim:Vbeta-Delta-2}
   \left\|\hat{V}_{\beta_2}(n,t,t_{0})-\hat{V}_{\beta_{1}}(n,t,t_{0})
   \right\|_n 
\ \le\ C(\eta,T) w(n)\int_{t_{0}}^t
   |\beta_2(t')-\beta_{1}(t')|\,dt'
\end{equation}
uniformly in $n$, which clearly implies Eq.~(\ref{claim:Vbeta-Delta}).
To see Eq.~(\ref{claim:Vbeta-Delta-2}), we write
\begin{equation}\label{eq:Vbeta-4}
\hat{V}_{\beta_2}(n,t,t_{0})-\hat{V}_{\beta_{1}}(n,t,t_{0})
\ =\ \int_{t_{0}}^t\hat{V}_{\beta_2}(n,t,t')\Delta \hat G(n,t')
\hat{V}_{\beta_{1}}(n,t',t_{0})dt'\ ,
\end{equation}
where 
\begin{equation}
\Delta \hat G(n,t)
\ =\ \left(\!\!\begin{array}{cc} 0& i(n\!+\!\beta_2(t))^3-i(n\!+\!\beta_{1}(t))^3\\
                    i(\beta_2(t)-\beta_{1}(t))&0 \end{array}\!\!\right)\ .
\end{equation} 
Using that $\abs{\beta_i(t)}\le 1+\eta T$ by assumption, 
we estimate
\begin{eqnarray}
\norm{\Delta \hat G(n,t)}_n 
&\le& \abs{\beta_2(t)-\beta_1(t)}\norm{
   \left(\!\!\begin{array}{cc}
  0& 3(|n|+1+\eta T)^2\\
   1&0 \end{array}\!\!\right)}_n\nonumber\\
&\le& \max{\left\{ 3\frac{(|n|+1+\eta T)^2}{w(n)}, w(n)\right\} }
         \abs{\beta_2(t)-\beta_1(t)} \nonumber \\
&\le& 3 (1+\eta T)^2 w(n) \abs{\beta_2(t)-\beta_1(t)}\ . 
        \label{eq:Vbeta-Delta-1}
\end{eqnarray}
Inserting Eq.~(\ref{eq:Vbeta-Delta-1}) and the bound 
in Eq.~(\ref{claim:Vbetabound-n}) of Lemma~\ref{lem:Vbetabound} 
into Eq.~(\ref{eq:Vbeta-4}) yields
Eq.~(\ref{claim:Vbeta-Delta-2}), with $C(\eta,T)=
c(\eta,T)^2 3 (1+\eta T)^2$, where $c(\eta,T)$ is the constant 
from Lemma~\ref{lem:Vbetabound}. This proves  the claim
in Eq.~(\ref{claim:Vbeta-Delta}).

The strong continuity follows
as in the proof of Lemma~\ref{lem:Vbetabound}
by combining Eq.~(\ref{claim:Vbeta-Delta-2}) with a uniform tail 
estimate obtained from Eq.~(\ref{claim:Vbetabound}).

\hfill $\Box$

\subsection{The resolvent for $-\partial_{s}^2 +\kappa^2$}

\bigskip\noindent
We provide a lower bound for the spectrum of the Schr\"{o}dinger
operator $L_{\kappa}\equiv-\partial_s^2 +\kappa^2$ acting in
$L^2[0,2\pi]$, with periodic boundary conditions.

\begin{lemma} \label{lem:spectrum}
The operator $L_{\kappa}$ has least eigenvalue 
$e_0(\kappa)$ satisfying,
\begin{equation}
e_0(\kappa)\geq 1/4\ .
\end{equation}
In particular, $L_\kappa$ is invertible in $L^2$, and 
$\norm{L_\kappa^{-1}}_{L^2}\le 4$. 
\end{lemma}

\begin{remark} {\em It is natural to conjecture that 
$\inf e_{0} (\kappa)$ is actually attained for $\x(s)$ 
a circle where $\kappa (s) =1$ and $e_0(\kappa)=1$.
An analogous result due to Harrell and Loss says that
the second-lowest eigenvalue of $-\partial_s^2-\kappa^2$
is maximal for a circle~\cite{HL}.}
\end{remark}

\medskip\noindent{\sc Proof of Lemma~\ref{lem:spectrum}:\ } \nobreak
We have that
\begin{eqnarray}\label{infeqapp}
 \inf_{\{\psi:\|\psi\|_2=1\}} (\psi, L_{\kappa}\psi)
&=&  \inf_{\{\psi:\|\psi\|_2=1\}}\left(\|\partial_s\psi\|_2^2
           +\|\psi \partial_s\u\|_2^2\right)\nonumber\\
&=& \inf_{\{\psi:\|\psi\|_2=1\}}\left(\|\partial_s(\psi\u)\|_2^2\right)\nonumber\\
&\ge& \inf_{\{\psivec\}} {' \|\partial_s \psivec\|_2^2}
\end{eqnarray}
where the vector function $\psivec= \psi\u$, and in the last
line, the infimum is to be taken over all normalized $\psivec$
such that
\begin{equation}
  \int \u(s) \, ds= \int \frac{\psivec(s)}{|\psivec(s)|} \,ds =
      0.
\end{equation}
If these constraint integrals are zero, then each component
${\psi_i}$ of $\psivec$ must vanish at some point $s_i$,
$i=1,2,3$, i.e., each component must satisfy a
Dirichlet condition and so $\|\partial_s\psi_{i}|_2^{2}\geq\frac{1}{4} \|\psi_{i}\|_2^{2}$, for each $i$ so that the infimum in the last line of
Eq.~(\ref{infeqapp}) is bounded below by
\begin{equation}
  \inf_{\sum \|\psi_i\|_2^2=1}{\sum \|\partial_s\psi_i\|_2^2}
\ \ge\ \inf_{\sum \|\psi_i\|_2^2=1}\frac{1}{4}\sum \|\psi_i\|_2^2
\ =\ \frac{1}{4}.
\end{equation}
This completes the proof of the eigenvalue estimate.

\hfill$\Box$

\bigskip
Next, we provide bounds on the resolvent for 
$L_{\kappa}$, considered as mapping $H^r$ to $H^{r+2}$,
assuming the spectral bound in Lemma~\ref{lem:spectrum}
and additional bounds on the norm of $\kappa$.
The inequalities given are by no means
optimal, but they are adequate for our purposes; the estimates are in
the spirit of Bessel kernel estimates (see \cite{Stein}). 

\begin{lemma} \label{lem:mu}
If $\kappa\in H^{r+1}$ for some $r\ge 0$, then 
$L_\kappa^{-1}$ defines a bounded linear operator from $H^{r-1}$ 
to $H^{r+1}$. More precisely, there exists a 
constant $C_1= C_1(r)$ such that if 
$\mu$ solves
\begin{equation}\label{appeq2}
     (-\partial_s^2+\kappa^2)\mu= f
\end{equation}
with $f\in  H^{r-1}$ for some $\kappa$ with $\norm{\kappa}_{H^{r+1}}\le R$,
then $\mu$  satisfies
\begin{equation}\label{appeq2.5}
\|\mu\|_{H^{r+1}} \ \leq \  C_1\left(1+
        R^2\right)^\nu \cdot \| f\|_{H^{r-1}}\ ,
\end{equation}
where $\nu=\nu(r)$ is the smallest the smallest integer 
at least as large as $(r+1)/2$ for $r\ge 1$,
and $\nu(r)=2$ for $0\le r<1$.  The Fourier coefficients of 
$\mu$ are bounded by
\begin{equation}
\label{appeq2.55}
\abs{\hat\mu(n)}
\ \le\ C_2 w(n)^{-2}\left(
    \abs{\hat f(n)} + w(n)^{-1-r}(1+R^2)^\nu\cdot 
                  \norm{f}_{H^{{r}-1}}\right)\ .
\end{equation}
\end{lemma}
 
\proof We compare $L_\kappa$ with the operator 
$-\partial_s^2 + 1$, using two forms of the resolvent identity:
\begin{eqnarray}
L_{\kappa}^{-1}
\label{reseq1}
&=&(-\partial_s^2+1)^{-1} - (-\partial_s^2+1)^{-1}(\kappa^2-1)L_\kappa^{-1}\\
\label{reseq2}
&=&(-\partial_s^2+1)^{-1} - L_\kappa^{-1}(\kappa^2-1)(-\partial_s^2+1)^{-1}\ .
\end{eqnarray}
Clearly,
\begin{equation}
\label{eq:L1}
\norm{(-\partial_s^2+1)^{-1} f}_{H^{{r'}+1}}\ = \ \norm{f}_{H^{r'}-1}
\end{equation} 
for any $r'\in\R$.  Let $f\in L^2$.
The first line of the resolvent identity shows that
\begin{eqnarray}
\nonumber
\norm{L_\kappa^{-1}f}_{H^2}
&\le& 
\norm{f}_{L^2} + \norm{(\kappa^2-1)L_\kappa^{-1}f}_{L^2}\\
&\le& c_1(1+R^2)\norm{f}_{L^2}\ 
\end{eqnarray}
for some constant $c_1$. We have used Eq.~(\ref{eq:L1}) in the first line.
In the second line, we have used Lemma~\ref{lem:Leibnitz}, which is proved
below, and the spectral bound of Lemma~\ref{lem:spectrum}. This proves 
the claim in the case
$r=1$.  More generally, the first line of the resolvent identity
shows that for $r'\le r$,
\begin{equation}
\norm{L_\kappa^{-1}}_{H^{r'+1}\to H^{{r'}+3}}\ \le\ 
c_2 (1+R^2)\norm{L_\kappa^{-1}}_{H^{r'-1}\to H^{r'+1}}
\end{equation}
for some constant $c_2=c_2(r,r')$.  
Iterating this estimate we obtain the claim for $r\ge 1$.
The case $0\le r\le 1$ follows from
the second line of the resolvent identity,
which implies that for $r'\le r$
\begin{equation}
\norm{L_\kappa^{-1}}_{H^{r'-3}\to H^{{r'}-1}}\ \le\ 
c_3 (1+R^2)\norm{L_\kappa^{-1}}_{H^{r'-1}\to H^{{r'}+1}}\ .
\end{equation}

\hfill $\Box$

We will also need estimates on how $\mu$ varies as the curvature 
$\kappa$ varies. 

\begin{lemma} \label{lem:mu-Delta}
Suppose that $\mu_1$ and $\mu_2$ solve the equations
\begin{eqnarray}
   L_{\kappa_1}\mu_1&=& (-\partial_s^2+\kappa_1^2)\mu_1= f\nonumber\\
   L_{\kappa_2}\mu_2&=& (-\partial_s^2+\kappa_2^2)\mu_2= f\ ,
\end{eqnarray}
where $\kappa_1,\kappa_2\in H^{r+1}$.
Suppose that 
$R\ge\max{\{\norm{\kappa_1}_{H^{r+1}},\norm{\kappa_2}_{H^{r+1}}\}}\,$.
For $r\ge 0$ there exists a constant $C=C(r)$ such that
\begin{equation}
\label{claim:mu-Delta} \norm{\mu_2-\mu_1}_{H^{r+1}}
\ \le\ C R \left(1+R^2\right)^{2\nu}\norm{\kappa_2-\kappa_1}_{H^r} 
                          \norm{f}_{H^{r-1}}\ .
\end{equation}
Here, $\nu=\nu(r)$ is the exponent from  Lemma~\ref{lem:mu}.
\end{lemma}

{F}rom the resolvent identity
\begin{equation}
\mu_2-\mu_1 
\ =\  -L_{\kappa_2}^{-1} (\kappa_2^{2}-\kappa_1^{2})L_{\kappa_1}^{-1}f\ ,
\end{equation}
we obtain for $r\ge 0$ with the help of Lemma~\ref{lem:Leibnitz},
which is proved below,
\begin{eqnarray}
\norm{\mu_2-\mu_1}_{r+1}
&\le& 2R \norm{L_{\kappa_2}^{-1}}_{H^{r}\to H^{r+1}}
         \norm{\kappa_2-\kappa_1}_{H^r}
         \norm{L_{\kappa_1}^{-1}}_{H^{r-1}\to H^r}
         \norm{f}_{H^{r-1}}\ .\quad
\end{eqnarray}
The claim now follows from Lemma~\ref{lem:mu}.

\hfill $\Box$

\subsection{The nonlinearity}

\begin{lemma} \label{lem:FB}
For $Y=(P,Q)\in \Y^r$ and $\beta \in \R$, let $F_\beta$ be 
defined by Eq.~(\ref{eq:def-F}), and let $B$ be defined by 
Eq.~(\ref{eq:def-B}).  There exist increasing continuous functions 
$C_1$, $C_2$, and $C_3$ of two variables,
which also depend on $r$, such that
the following estimates hold:
\begin{enumerate} \item
If $r\ge 0$, then for any $Y\in\Y^r$ with
$\norm{Y}_{\Y^r}\le R$, and any $\beta\in\R$ with $\abs{\beta}\le b$,
\begin{equation} \label{claim:FB-bound}
\left\{ \!\!
\begin{array}{rcl}
\norm{F_\beta(Y)}_{\Y^r} &\le& C_1(b,R) \\
\,\abs{B(Y)}               &\le& R^2\ .
\end{array}\right.
\end{equation}
\item If $r\ge 0$, then for any pair 
$\Y_1,Y_2\in \Y^r$ with $\norm{Y_i}_{\Y^r}\le R$
and any $\beta_1,\beta_2\in\R$ with $\abs{\beta_i}\le b$,
\begin{equation} \label{claim:FB-Lip}
\left\{ \!\!
\begin{array}{rcl}
\norm{F_{\beta_2}(Y_2) \!- \!F_{\beta_1}(Y_1)}_{\Y^r}
                    &\le& C_2(b,R) \Bigl( \norm{Y_2\!-\!Y_1}_{\Y^r}
                            + \abs{\beta_2\!-\!\beta_1}\Bigr)\\
\,\abs{B(Y_2)\!-\!B(Y_1)} &\le& R\|Y_2\!-\!Y_1\|_{\Y^r}\ .
\end{array} \right.
\end{equation}
\item If $r\ge 1/2$, then for any pair $Y_1,Y_2\in \Y^r$ 
with $\norm{Y_i}_{\Y^r}\le R$
and any $\beta_1,\beta_2\in\R$ with $\abs{\beta_i}\le b$,
\begin{equation} \label{claim:FB-Delta}
\left\{ \!\!
\begin{array}{rcl}
\norm{F_{\beta_2}(Y_2) \!- \!F_{\beta_1}(Y_1)}_{\Y^{r-1}}\!\!
   &\le&C_3(b,R) \Bigl(\norm{Y_2\!-\!Y_1}_{\Y^{r\!-\!1}} 
                    + \abs{\beta_2\!-\!\beta_1}\Bigr)\\
\,\abs{B(Y_2)\!-\!B(Y_1)} &\le& R\|Y_2\!-\!Y_1\|_{\Y^{r\!-\!1}}\ .
\end{array}\right.
\end{equation}
\end{enumerate} 
\end{lemma}

\noindent The following lemma will be useful in estimating
the various terms in the nonlinearity.

\begin{lemma} {\rm (Leibnitz rule for $H^r$-norms)} \ 
\label{lem:Leibnitz}
There exist constants $C_1$ - $C_5$ which
depend only on $r$, so that each of the following inequalities
holds whenever the right hand side is finite:
\begin{equation}
\begin{array}{lcll}
\norm{fg}_{H^r}     &\le& C_1\norm{f}_{H^r}\norm{g}_{H^r}\quad &(r\ge 1)\\
\norm{fg}_{H^r}     &\le& C_2\norm{f}_{H^r}\norm{g}_{H^{r+1}}\quad& (r\ge 0)\\
\norm{fg}_{H^{r-1}} &\le& C_3\norm{f}_{H^{r}}\norm{g}_{H^{r}}\quad&(r\ge 0)\\
\norm{fg}_{H^{r-1}} &\le& C_4\norm{f}_{H^{r-1}} \norm{g}_{H^{r+1}}\quad 
                                         &(r\ge 1/2)\\
\norm{fg}_{H^{r-2}} &\le& C_5\norm{f}_{H^{r-1}} \norm{g}_{H^{r}} \quad 
                                         &(r\ge 1/2)\ .
\end{array}
\label{claim:Leibnitz}
\end{equation}
\end{lemma}

\proof If $r$ is an integer, then 
the first two inequalities 
follow immediately from the Leibnitz rule and
the fact that $H^1\subset L^\infty$ in one space dimension.  
In general, we use that, for $r\ge 0$,
$w(n)^r\le c_1(w(k)^r + w(n-k)^r)$
with some constant $c_1=c_1(r)$ to get
\begin{eqnarray}
w(n)^r \abs{\widehat{fg}(n)}
&\le& w(n)^r\sum_k \abs{\hat f(k)}\cdot\abs{\hat g(n-k)}\\
&& \hskip -3cm  c_1 \left(\!
   \sum_k w(k)^r \abs{\hat f(k)} \cdot \abs{\hat g(n-k)}
   + \sum_k\abs{\hat f(k)} \cdot 
         w(n-k)^r\abs{\hat g(n-k)}\!\right)\,.
\nonumber
\end{eqnarray}
If either $f,g\in H^r$ with $r\ge 1$, or $f\in H^r,g\in H^{r+1}$,
each of the two sums on the right hand side is the Fourier
transform of the product of an $L^2$-function
with an $H^1$-function, and hence in $\ell^2$.
This proves the first two inequalities.

For $r\ge 1$, 
the third inequality follows from the first.
For $0\le r<1$, we use that 
$w(n)\le \sqrt{2} w(k)w(n-k)$ to obtain
\begin{equation}
w(n)^{r-1}\abs{\widehat{fg}(n)}
\ \le\ c_2 w(n)^{-1} \sum_k w(k)^r\abs{\hat f(k)}
         \cdot w(n-k)^r\abs{ \hat g(n-k)} \ ,
\end{equation}
where the constant $c_2$ depends only on $r$.
The sum is the Fourier transform of the product
of two  $L^2$-functions, and hence 
in $\ell^\infty$. Since $w^{-1}\in\ell^2$,
it follows that $w(n)^{r-1}\widehat{fg}\in \ell^2$, which implies the 
third inequality.

For $r\ge 1$, the last two inequalities
follow from the second and third with $r$ replaced by $r-1$.  
For $1/2\le r < 1$, we use that $\sqrt{2}w(n)\ge w(k)/w(n-k)$ and 
proceed as in the proof of the third inequality.

\hfill$\Box$

\bigskip\noindent{\sc Proof of Lemma~\ref{lem:FB}} \ 
The estimates for $F^{(1)}$, $F^{(3)}$, and $B$ follow
by repeated applications of Lemma~\ref{lem:Leibnitz}. We focus
on the terms involving $\mu$.

Consider the first claim, Eq.~(\ref{claim:FB-bound})
for $r\ge 0$. We want to apply 
Lemma~\ref{lem:mu} with $\kappa=\abs{Q}\in H^{r+1}$
and $f= \abs{P}^2 + \abs{Q}^4 - \abs{(\partial_s + i\beta)Q}^2$.
By the third inequality of Lemma~\ref{lem:Leibnitz}, we have for $r\ge 0$,
\begin{equation}
\|f\|_{H^{r-1}}\ \le\ c_1 (1+\abs{\beta}^2)R^2\ .
\end{equation}
By Lemma~\ref{lem:mu},
\begin{equation}
\norm{\mu}_{H^{r+1}}\ \le\ c_2 (1+\abs{\beta}^2) R^2 (1+R^2)^\nu\ ,
\end{equation}
and so by the second inequality of Lemma~\ref{lem:Leibnitz},
\begin{eqnarray}
\|F_\beta^{(2)}(Y)\|_{\Y^r}
&=& \nonumber \norm{2(\partial_s\mu)Q + \mu(\partial_s+i\beta)Q}_{H^r}\\
&\le& c_3(1+\abs{\beta}) \norm{\mu}_{H^{r+1}}\norm{Q}_{H^{r+1}}\\
&\le& \nonumber c_4(1+\abs{\beta}^3) R^3 (1+R^2)^\nu\ ,
\label{eq:F-mu}
\end{eqnarray}
where $\nu=\nu(r)$ is the exponent from Lemma~\ref{lem:mu}.
This shows the bound in Eq.~(\ref{claim:FB-bound}).
The proof of the second claim, Eq.~(\ref{claim:FB-Lip}), 
is almost the same.

To see the third claim, Eq.~(\ref{claim:FB-Delta}), let
$Y_i=(P_i,Q_i)$, $\kappa_i=\abs{Q_i}$, and denote 
by $f_i$ the right hand side of Eq.~(\ref{eq:def-mu}) corresponding 
to $Y_i$ ($i=1,2$).  By Lemma~\ref{lem:Leibnitz}, we have for $r\ge 1/2$,
\begin{eqnarray}
\norm{f_2\!-\!f_1}_{H^{r\!-\!2}}  \!
&\le& c_1 (1+b^2)R\norm{Y_1\!-\!Y_2}_{\Y^{r\!-\!1}} 
             + c_2 b R^2 \abs{\beta_2\!-\!\beta_1}\\
&\le& \nonumber c_3(1+b^2)(1+R^2) \left(\norm{Y_2\!-\!Y_1}_{\Y^{r\!-\!1}} + 
                \abs{\beta_2\!-\!\beta_1}\right)
\end{eqnarray}
with suitable constants $c_1$-$c_3$
By Lemmas~\ref{lem:mu} and~\ref{lem:mu-Delta}, this implies
\begin{eqnarray}
\nonumber 
\norm{\mu_2-\mu_1}_{H^{r}} 
&\le& \norm{L_{\kappa_2}^{-1}(f_2-f_1)}_{H^{r}} + 
       \norm{\left(L_{\kappa_2}^{-1}- L_{\kappa_1}^{-1}\right) f_1}_{H^r}\\
&\le& \nonumber c_4(1+b^2)(1+R^2)^\nu \norm{f_2-f_1}_{H^{r-2}}\\
&& \quad + 
     c_5R(1+R^2)^{2\nu}\norm{\kappa_2-\kappa_1}_{H^r}\norm{f}_{H^{r-1}}\\
&\le& \nonumber c_6(1+b^2)(1+R^2)^{2\nu+3/2}\Bigl(\norm{Y_2-Y_1}_{\Y{r-1}}
          +    \abs{\beta_2-\beta_1}\Bigr)
\end{eqnarray}
with suitable constants $c_4$-$c_6$.
Inserting this and Eq.~(\ref{eq:F-mu})
into $F_\beta^{(2)}$ and using the fourth inequality
of Lemma~\ref{lem:Leibnitz}, we see that 
$F^{(2)}$ satisfies the bound in Eq.~(\ref{claim:FB-Delta}).
\hfill$\Box$

\section{Proof of Theorem~\ref{thm:main}}
\label{sect:local}
\setcounter{equation}{0}

\subsection{Existence of solutions to 
Eqs.~(\ref{eq:PQ-all})-(\ref{eq:initial-PQbeta})}

In the notation of Eqs.~(\ref{eq:def-G})-(\ref{eq:PQ-linear}),
the Duhamel formula for Eqs.~(\ref{eq:PQ-all})-(\ref{eq:initial-PQbeta})
is given by
\begin{equation} \label{eq:Y-Duhamel}
\left\{
\begin{array}{rclcl}
Y(t) &=& V_\beta(t,0) Y_0
           + \int_0^t V_\beta(t,t') F_{\beta(t')}(Y(t')) \,dt' 
       &=:& {\cal F}_\beta(Y)(t)\\[0.1cm]
\beta(t) &=& \beta_0  +  \int_0^t B(Y(t'))\, dt' 
       &=:& {\cal B}(Y)(t)\ ,
\end{array}\right.
\end{equation}
where $Y_0=(P_0,Q_0)\in \Y^r$ and $\beta_0\in \R$ are given initial values.
We begin by solving the first equation in Eq.~(\ref{eq:Y-Duhamel})
for a fixed function $\beta$.

\begin{lemma} \ \label{lem:fixed-beta}
Let $\beta$ be a Lipschitz continuous real-valued function 
on $\R^+$ with Lipschitz constant $\eta$ and  $\abs{\beta(0)}\le 1$, and
let $Y_0\in\Y^r$ for some $r\ge 0$.
There exists a number $R<\infty $, which depends on 
$\|Y_0\|_{\Y^r}$, and a time $T>0$ which depends on $\eta$ and $R$
such that the fixed point equation
\begin{equation}
\label{eq:fix-fixed-beta}
Y={\cal F}_\beta(Y)
\end{equation}
has a unique solution $Y^*_\beta$ on $[0,T]$  in $\Y^r$.
The solution  satisfies the uniform bound
\begin{equation}
\sup_{0\le t\le T} \norm{Y^*_\beta(t)}_{\Y^r}\ \le R\ ,
\end{equation}
and depends continuously on the initial value $Y_0$ with respect to the
$\Y^r$-norm.  It is also strongly continuous in $\beta$ with respect
to the $\Y^r$-norm, uniformly on $[0,T]$, in the sense
that for any $\eps$ there exists $\delta=\delta(\eps,\beta, T,\Y_0)$ 
such that
\begin{equation}
\sup_{0\le t\le T}\abs{\beta'(t)-\beta(t)}\  \le\ \delta
\end{equation}
implies that the corresponding solutions $Y^*_{\beta}$ and 
$Y^*_{\beta'}$ satisfy
\begin{equation}
\sup_{0\le t\le T}\|Y^*_{\beta'}(t)-Y^*_{\beta}(t)\|_{\Y^r}\ \le\ \eps\ .
\end{equation}
\end{lemma}

\proof We consider ${\cal F}_\beta$ as a map from the space
\begin{equation}
{\cal D}_R\ =\ 
\left\{ Y:[0,T]\mapsto \Y^r\ \big\vert \ Y\ \mbox{continuous},
\ \sup_{0\le t\le T}\|Y(t)\|_{\Y^r}\le R \right\}\ 
\end{equation}
with norm
\begin{equation}
|||Y|||_{r}\ =\ \sup_{0\le t\le T} \|Y(t)\|_{\Y^r}
\end{equation}
into the space of real-valued continuous functions on $[0,T]$ with
values in $\Y^r$.  The values of $T$ and $R$ will be chosen below.

By Lemma~\ref{lem:Vbetabound}, there exists an increasing continuous 
function of two variables $C_1$  with $\sup_\eta C_1(\eta,0)<\infty$ 
such that for any $Z\in \Y^r$, we have
\begin{equation}
\sup_{0\le t_0 \le  t\le T}\norm{V_\beta(t,t_0)Z}_{\Y^r}\ 
       \le C_1(\eta,T) \norm{Z}_{\Y^r}\ .
\end{equation}
By Eq.~(\ref{claim:FB-Lip}) of Lemma~\ref{lem:FB}, there exists an increasing
continuous function $C_2$ of two variables such that for
any two continuous functions $Y_1,Y_2$ on $[0,T]$ 
with values in  $\Y^r$ which are bounded uniformly by $R$,
\begin{equation}
\norm{F_{\beta(t)}(Y_2(t))-F_{\beta(t)}(Y_1(t))}_{\Y^r}\ \le\
        C_2(\abs{\beta(t)},R) \norm{Y_2(t)-Y_1(t)}_{\Y^r}\ 
\end{equation}
for all $0\le t\le T$. Combining the above two estimates
and using that $\abs{\beta(t)}\le 1+\eta T$, we see that
\begin{eqnarray}
\nonumber 
\sup_{0\le t\le T}
           \norm{{\cal F}_\beta (Y_2)(t) -{\cal F}_\beta(Y_1)(t)}_{\Y^r}
\hskip -3cm &&\\
&\le& \int_0^T C_1(\eta,T)
      \norm{F(Y_2(t')) - F(Y_1(t'))}_{\Y^r} \, dt' \\
 &\le& T C_1(\eta,T) C_2(1+\eta T,R) 
        \sup_{0\le t\le T} \norm{Y_2(t)-Y_1(t)}_{\Y^r}\ .
\nonumber
\end{eqnarray}
Note that by Lemma~\ref{lem:Vbetabound} and Eq.~(\ref{claim:FB-Lip})
of Lemma~\ref{lem:FB}, the function ${\cal F}_\beta(Y)$ is again
a continuous function of $t$. 

Fix
\begin{equation}
\label{eq:choose-R}
R\ \ge\  4 \sup_{\eta\in\R}C_1(\eta,0)\norm{Y_0}_{\Y^r}\ ,
\end{equation}
and choose $T$ small enough so that
\begin{equation}
\label{eq:choose-T}
\quad C_1(\eta,T)\ \le\ 2C_1(\eta,0)\ ,
\quad TC_1(\eta,T)C_2(1+\eta T,R) \ \le\ \frac{1}{2} \ .
\end{equation}
Then ${\cal F}_\beta$ is a contraction with Lipschitz constant 
$1/2$ on ${\cal D}_R$. Since, for \mbox{$Y\in {\cal D}_R$,}
\begin{eqnarray}
\sup_{0\le t\le T} \|{\cal F}_\beta(Y)(t)\|_{\Y^r}
\!&\le&\! \sup_{0\le t\le T} \|V_\beta(t,0)Y_0\|_{\Y^r}
       + \sup_{0\le t\le T}\|{\cal F}_\beta(Y)\!-\!{\cal F}_\beta(0)\|_{\Y^r}
\nonumber \\
&\le&\! R \ ,
\end{eqnarray}
we see that ${\cal F}_\beta$ maps ${\cal D}_R$ into itself.
By the contraction mapping principle, ${\cal F}_\beta$ has a unique 
fixed point in ${\cal D}_R$, which we denote by $Y^*_\beta$. 

The function ${\cal F}_\beta$ is clearly continuous in the 
initial value $Y_0$ with respect to the norm on ${\cal D}_R$. 
By Lemma~\ref{lem:Vbeta-Delta} and 
Eq.~(\ref{claim:FB-Lip}) of Lemma~\ref{lem:FB}, it
is also strongly continuous in $\beta$, in the sense that
for every $\eps>0$ there exists $\delta=\delta(\eps,\eta,T,Y)$
such that
\begin{equation}
\sup_{0\le t\le T} \abs{\beta'(t)-\beta(t)}\ \le\ \delta
\end{equation}
implies
\begin{equation}
   \sup_{0\le t\le T} \|{\cal F}_{\beta'}(Y) - {\cal F}_{\beta}(Y)\|_{\Y^r}
    \ \le\ \eps\ ,
\end{equation}
where $Y$ is any continuous function on $[0,T]$ with
values in $\Y^r$, and  $\beta$ is Lipschitz continuous with
$\abs{\beta(0)}\le 1$ and Lipschitz constant $\eta$.
By the uniform contraction principle (see~\cite{Chow-Hale}, Theorem 2.2), 
the fixed point $Y^*_{\beta}$ inherits the claimed 
continuity  properties from ${\cal F}_\beta$.
The modulus of continuity depends
on $\beta$ and on $Y_0$ through  the dependence
of $Y^*_{\beta}$ on these parameters.

\hfill$\Box$

\begin{coro} \label{cor:planar}
Let $Y_0=(P_0,Q_0)\in \Y^r$ for some $r\ge 0$,
where $P$ and $Q$ are real-valued, and let $\beta_0=0$.
There exists a time $T>0$, which depends on $\|Y_0\|_{\Y^r}$,
such that the fixed point equations in Eq.~(\ref{eq:Y-Duhamel})
have a unique solution of the form $(Y^*,0)$ on $[0, T]$,
where $Y^*=(P^*,Q^*)\in \Y^r$ are real-valued.
The solution depends continuously on the initial values in $\Y^r$.
\end{coro}

\proof The claim follows immediately from 
Lemma~\ref{lem:fixed-beta} with $\beta\equiv 0$.

\begin{theorem} \label{thm:exist}
Let $r\ge 0$, and let initial values $Y_0\in\Y^r$ and $\beta_0\in\R$
be given. There exists a time $T>0$,
which depends on $\norm{Y_0}_{\Y^r}$ and $\abs{\beta_0}$,
such that the pair of fixed point equations 
in Eq.~(\ref{eq:Y-Duhamel}) has a solution in $\Y^r\times\R$ 
on $[0,T]$ which is strongly continuous in $t$ and 
assumes the initial values at $t=0$.
\end{theorem}

\proof We may assume by the gauge invariance
given in Eq.~(\ref{eq:PQ-gauge}) that $\abs{\beta_0} < 1$.
Let $R=R(\|Y_0\|_{\Y^r})$ be the constant appearing in the 
statement of Lemma~(\ref{lem:fixed-beta}), and consider the space
\begin{equation}
{\cal C}_R\ =\ \{ \beta:[0,T]\mapsto \R\ \mid\ \abs{\beta(0)}\le 1,\ 
          \mbox{Lip}(\beta)\le R^2 \}\ 
\end{equation}
with the norm of uniform convergence. Then ${\cal C}_R$ is a compact 
convex subset of the space of continuous real-valued functions 
on $[0,T]$.

By  Lemma~\ref{lem:fixed-beta}, there exists a time $T=T(R)>0$
such that for $\beta\in {\cal C}_R$, the fixed point equation
in Eqs.~(\ref{eq:fix-fixed-beta}) has a unique solution $Y^*_\beta$
which is continuous in $t$ and satisfies
\begin{equation}
\sup_{0\le t\le T} \|Y^*_\beta(t)\|_{\Y^r}\ \le\ R\ .
\end{equation}
By Eq.~(\ref{claim:FB-bound}) of Lemma~(\ref{lem:FB}), we have that 
\begin{equation}
\abs{ {\cal B}(Y^*_\beta)(t)- {\cal B}(Y^*_\beta)(t_0)}
\ \le\ \int_{t_0}^t B(Y^*_\beta)(t')\, dt'
\ \le\ R^2\abs{t-t_0}\ ,
\end{equation}
that is, $B(Y^*_\beta)$ is again Lipschitz continuous with
Lipschitz constant $R^2$, and hence in ${\cal C}_R$.
Since the map $\beta\mapsto {\cal B}(Y^*_\beta)$
is continuous on ${\cal D}_R$ by Lemma~\ref{lem:fixed-beta}
and Eq.~(\ref{claim:FB-Lip}) of Lemma~\ref{lem:FB}, Schauder's theorem
implies that it has a fixed point in ${\cal C}_R$, which we denote by
$\beta^*$.  By construction, the pair 
$\left(Y^*_{\beta^*},\beta^*\right)$ solves the fixed point 
equations in Eq.~(\ref{eq:Y-Duhamel}).

\hfill$\Box$

\begin{remark}  {\em The dependence of $R$ and $T$ on the
initial value $\beta_0$ is due to the fact that
the gauge transformation  in Eq.~(\ref{eq:PQ-gauge})
may change  the norm of the initial values $Y_0$.
Let $k$ be the integer closest to $\beta_0$, and perform the gauge 
transformation in Eq.~(\ref{eq:PQ-gauge}) with $s_0=0$. Then
$\abs{\beta_0-k}\le 1/2$, and for $r\ge 0$,
\begin{equation}
\norm{e^{iks}Y_0}_{\Y^r}\ \le\ (1+k^2)^{(r+1)/2}\norm{Y_0}_{\Y^r}\ .
\end{equation}
}
\end{remark}

\subsection{Well-posedness}

In order to exploit the equivalence of 
Eqs.~(\ref{eq:PQ-all})-(\ref{eq:initial-PQbeta})
with the original initial value problem 
in Eqs.~(\ref{eq:EL-u})-(\ref{eq:initial-u})
that was established in Lemma~\ref{lem:equiv},
we need to prove that the solutions
of Eqs.~(\ref{eq:PQ-all})-(\ref{eq:initial-PQbeta})
are unique, and depend continuously on the initial values 
in some reasonable norm. The difficulty is that the linear 
operator $V_\beta$ is only strongly continuous, not norm-continuous,
with respect to the parameter $\beta$ (see 
Lemma~\ref{lem:Vbeta-Delta}).  We deal with this problem by 
using a weaker norm to estimate ${\cal F}_\beta$ and ${\cal B}$.

\begin{theorem} \label{thm:unique}
Let $r\ge 1/2$. For given initial values 
$Y_0\in\Y^r$ and $\beta_0\in\R$, there exists a time $T>0$ which
depends on $\|Y_0\|_{\Y^r}$ and on $\abs{\beta_0}$
such that the pair fixed point equations 
in Eq.~(\ref{eq:Y-Duhamel}) has a unique solution 
in $\Y^r\times \R$ on $[0,T]$ which is strongly continuous in $t$.
The solution depends continuously on the initial 
values with respect to the norm on $\Y^r\times\R$.
\end{theorem}

\proof By the gauge invariance in Eq.~(\ref{eq:PQ-gauge}),
we may assume that $\abs{\beta_0}<1$. 
For $R$ to be chosen below, we consider the right hand side 
of Eq.~(\ref{eq:Y-Duhamel}) as a map on the space
\begin{equation}
{\cal D}_R\!\times\! {\cal C}_R
\ =\ \left\{ (Y,\beta): [0,T]\mapsto \Y^r\ \Bigg\vert\ 
\begin{array}{cc}
Y\ \mbox{continuous},\\[0.1cm]
\sup_{0\le t\le T} \|Y(t)\|_{\Y^r}\le R,\\[0.1cm]
\mbox{Lip}(\beta)\le R^2,\ \abs{\beta(0)}\le 1 \end{array}
\right\}\ ,
\end{equation}
with the norm
\begin{equation}
|||(Y,\beta)|||\ =\ 
\sup_{0\le t\le T} \|Y(t)\|_{\Y^{r-1}}  + \sup_{0\le t\le T} \abs{\beta(t)}\ .
\end{equation}
Note that ${\cal D}_R\!\times\! {\cal C}_R$ is complete with respect
to the $|||\cdot|||$-norm by the convexity 
of the $\Y^r$-norm and the fact that the uniform limit of 
Lipschitz continuous functions with a given constant
is again Lipschitz continuous with that constant.

We first bound the $\Y^r$-norm
of ${\cal F}_\beta(Y)$ and the Lipschitz constant of ${\cal B}(Y)$.
By Lemma~\ref{lem:Vbetabound}, there exists an 
increasing continuous function $C_1$ of two variables with 
$\sup_{\eta} C_1(\eta,0)<\infty$ such that for any
Lipschitz continuous function $\beta$ and any $Z\in\Y^r$,
\begin{equation}
\sup_{0\le t_0\le t\le T}\norm{V_\beta(t,t_0)Z}_{Y^r} 
\ \le\ C_1(\mbox{Lip}(\beta),T)  \|Z\|_{\Y^r} 
\label{eq:bound-1}
\end{equation}
By Eq.~(\ref{claim:FB-bound}) of Lemma~\ref{lem:FB}, there exists 
an increasing continuous function $C_2$ of two variables such that
for any $(Y,\beta)\in {\cal D}_R\!\times\!{\cal C}_R$,
\begin{equation}
\sup_{0\le t\le T}\norm{F_{\beta(t)}(Y(t))}_{\Y^{r}}
       \ \le\ C_2(\abs{\beta(t)},R) \ .
\label{eq:bound-2}
\end{equation}
Combining Eqs.~(\ref{eq:bound-1})-(\ref{eq:bound-2}),
we see that for $(Y,\beta)\in {\cal D}_R\!\times\!{\cal C}_R$,
\begin{eqnarray} 
\label{eq:bound-3}
\sup_{0\le t\le T} \norm{{\cal F}_\beta(Y,\beta)(t)}_{\Y^r}
&\le & C_1(R^2,T) \norm{Y_0}_{\Y^r}\\
&&\quad  +  T C_1(R^2,T) C_2(1+R^2T,R)\ . \nonumber
\end{eqnarray}
By Eq.~(\ref{claim:FB-bound}) of Lemma~\ref{lem:FB}, we also have 
\begin{equation}
{\rm Lip}\bigl( {\cal B}(Y)\bigr)\ =\ \sup_{0\le t_0<t\le T}
\frac{\abs{ {\cal B}(Y)(t)-{\cal B}(Y)(t_0)}}{\abs{t-t_0}}\ \le\ R^2\ .
\label{eq:bound-4}
\end{equation}

Next we bound the Lipschitz constant of $({\cal F}_\beta,{\cal B})$
on ${\cal X}$.  
Let $(Y_1,\beta_1)$ and $(Y_2,\beta_2)$ be in ${\cal D}_R\times {\cal C}_R$.
By Lemma~\ref{lem:Vbeta-Delta}, there exists 
an increasing continuous function $C_3$
of two variables such that for any two real-valued 
functions $\beta_1$, $\beta_2$
which are Lipschitz continuous with Lipschitz constant $R^2$,
\begin{equation}
\norm{V_{\beta_2}(t,t_0)Z \!-\!V_{\beta_1}(t,t_0)Z}_{Y^{r\!-\!1}}
\ \le\ T C_3(R^2,T)\sup_{0\le t'\le T} 
           \abs{\beta_2(t')\!-\!\beta_1(t')} \norm{Z}_{\Y^r}
\label{eq:Delta-1}
\end{equation}
for all $0\le t_0\le t\le T$.
By Eq.~(\ref{claim:FB-Delta})
of Lemma~\ref{lem:FB}, there exists an increasing continuous 
function $C_4$ of two variables such that
\begin{eqnarray}
\nonumber
\norm{F_{\beta_2(t)}(Y_2(t))-F_{\beta_1(t)}(Y_1(t))}_{\Y^{r-1}}
&\le& C_4(\abs{\beta},R) 
  \Bigl( \norm{Y_2(t)-Y_1(t)}_{\Y^{r-1}} \\
&& \quad + \abs{\beta_2(t)-\beta_1(t)} \Bigr)
\label{eq:Delta-2}
\end{eqnarray}
for all $0\le t\le T$. 
Combining Eqs.~(\ref{eq:bound-1})-(\ref{eq:bound-2})
with Eqs.~(\ref{eq:Delta-1})-(\ref{eq:Delta-2}), and using that 
$\abs{\beta(t)}\le 1+R^2T$,
we see that for $(Y,\beta)\in {\cal D}_R\times{\cal C}_R$,
\begin{eqnarray}
\sup_{0\le t\le T}\norm{{\cal F}_{\beta_2}(Y_2,\beta_2)(t) -
       {\cal F}_{\beta_1}(Y_1,\beta_1)(t)}_{\Y^{r-1}} &&
\label{eq:Delta-3} \\
&& \hskip -4cm \le\ TC_5(R,T, \norm{Y_0}_{\Y^r}) 
|||(Y_2,\beta_2)-(Y_1,\beta_1)||| \nonumber
\end{eqnarray}
with some increasing continuous function $C_5$ of two variables. 
By Eq.~(\ref{claim:FB-Delta}) of Lemma~\ref{lem:FB},
\begin{equation}
\sup_{0\le t\le T}\abs{{\cal B}(Y_2)(t)-{\cal B}(Y_1)(t)}
\ \le\ TR \norm{Y_2-Y_1}_{\Y^{r-1}}\ .
\label{eq:Delta-4}
\end{equation}
In summary, Eqs.~(\ref{eq:Delta-3})-(\ref{eq:Delta-4})
show that
\begin{eqnarray}
\label{eq:Delta-5}
|||({\cal F}_{\beta_2}(Y_2),{\cal B}(Y_2)) 
-({\cal F}_{\beta_1}(Y_1),{\cal B}(Y_1))||| &&\\
&& \hskip -3cm \le\ TC_6(R,T,\|Y_0\|_{\Y^r})|||(Y_2,\beta_2)-(Y_1\beta_1)|||
\nonumber 
\end{eqnarray}
with a suitable increasing function $C_6$.

Choose 
\begin{equation}
\label{eq:choose-R-2}
R\ \ge\  4 \sup_{\eta}C_1(\eta,0) \norm{Y_0}_{\Y^r} \ ,
\end{equation}
and $T$ small enough such that
\begin{equation}
\label{eq:choose-T-2}
\begin{array}{c}
C_1(R^2,T)\ \le\ 2C_1(R^2,0)\ ,\\
TC_1(R^2,T)C_2(1+R^2T,R)\ \le\ R\ , \\
TC_6(R,T) \ \le \  \frac{1}{2} \ .
\end{array}
\end{equation}
Then Eqs.~(\ref{eq:bound-3})-(\ref{eq:bound-4})
show that $\bigl({\cal F}_\beta,{\cal B}\bigr)$ maps 
${\cal D}_R\times {\cal C}_R$ 
into itself, and Eq.~(\ref{eq:Delta-5})
shows that $({\cal F}_\beta,{\cal B})$ is
a contraction with Lipschitz constant $1/2$.
By the contraction mapping theorem, Eq.~(\ref{eq:Y-Duhamel}) has 
a unique solution in ${\cal D}_R\times{\cal C}_R\,$.  
By the uniform contraction principle 
and the continuity properties of $V_\beta(t,t')$ and $F_\beta$,
this solution is continuous with respect 
to the initial data $(Y_0,\beta_0)$ in 
the sense that for every $\eps>0$ there exists
$\delta=\delta(\eps,T,\beta_0,Y_0)$ such that
\begin{equation}
\abs{\beta_0'-\beta_0} + \|Y_0'-Y_0\|_{\Y^{r-1}}\ \le\ \delta
\end{equation}
implies that the corresponding solutions 
$(Y,\beta)$ and $(Y',\beta')$ of Eq.~(\ref{eq:Y-Duhamel})
satisfy
\begin{equation}
|||Y'-Y|||\ \le\ \eps\ .
\label{eq:r-1-continuous}
\end{equation}

It remains to prove the strong  continuity of the solution with 
respect to $Y_0$ and $\beta_0$ in the natural norm on $\Y^r\times \R$. 
Note that Eq.~(\ref{eq:r-1-continuous}) implies in particular that
\begin{equation}
\sup_{0\le t\le T} \abs{\beta'(t)-\beta(t)}\ \le\ \eps\ .
\label{eq:beta-continuous}
\end{equation}
Fix initial values $(Y_0,\beta_0)$ in $\Y^r\times\R$
with $\abs{\beta_0}<1$,
let $\eps>0$ be given, and suppose that
\begin{equation}
\abs{\beta_0'-\beta_0} + \|Y_0'-Y_0\|_{\Y^{r}}\ \le\ \delta
\end{equation}
Let $(Y, \beta)\in {\cal D}_R\times{\cal C}_R$ be the solution
of Eq.~(\ref{eq:Y-Duhamel}) with these initial values,
defined on some interval strictly containing $[0,T]$,
and let $(Y',\beta')$ be the solution to 
Eq.~(\ref{eq:Y-Duhamel}) with initial values $(Y'_0,\beta'_0)$,
where $Y_0\in \Y^r$, and $\abs{\beta_0}\le 1$. 
If $Y_0'$ is sufficiently close to $Y_0$,
then we may assume that $(Y'\beta')$ is defined on $[0,T]$. 
By Lemma~\ref{lem:fixed-beta}, $Y$ is the unique solution of
the first fixed point equation in Eq.~(\ref{eq:Y-Duhamel}) with 
$\beta$ fixed and the given initial value $Y_0$.
Let $Z$ be the unique solution of the fixed point equation
\begin{equation}
Z(t)\ =\ V_{\beta}(t,0)Y_0' + \int_0^t V_\beta(t,t')F_{\beta(t')}(Z(t'))\, dt'\ 
\end{equation}
with the same function $\beta$. Clearly,
\begin{equation}
\sup_{0\le t\le T}\|Y'(t)\!-\!Y(t)\|_{\Y^r}
\ \le\ \sup_{0\le t\le T}
   \Bigl(\|Y'(t)\!-\!Z(t)\|_{\Y^r} + \|Z(t)\!-\!Y(t)\|_{\Y^r}\Bigr)\,.
\end{equation}
By the continuity statement in Eq.~(\ref{eq:beta-continuous}),
we may assume that $\sup_{0\le t\le T} \abs{\beta'(t)-\beta(t)}$
is as small as we please.
Since $Y'$ and $Z$ solve the fixed point equation in the
first line of Eq.~(\ref{eq:Y-Duhamel}) with the same initial 
value $Y'_0$ but different functions $\beta'$ and $\beta$, and
since $Z$ and $Y$ solve the equation with the same function $\beta$ but 
different initial values $Y_0'$ and $Y_0$, the continuity statements
of Lemma~\ref{lem:fixed-beta} imply the claim.
The modulus of continuity depends on $\beta_0$ and $Y_0$ through
the dependence of $(Y,\beta)$ on these parameters.

\hfill$\Box$

\noindent{\sc Proof of Theorem~\ref{thm:main}} \ 
By Lemma~\ref{lem:equiv}, the well-posedness 
of Eqs.~(\ref{eq:EL-u})-(\ref{eq:initial-u}) asserted
in Theorem~\ref{thm:unique} for the three-dimensional case
and in Corollary~\ref{cor:planar} for the planar case
implies the claims of Theorem~\ref{thm:main}.

\hfill $\Box$

\bigskip \noindent {\bf Acknowledgments:} The authors wish to thank
J.G.Simmonds for his helpful remarks and for providing us with 
reference information. A.B. wishes to thank J. Colliander 
for a stimulating conversation.  L.E.T. wishes to acknowledge the
hospitality of the \'{E}cole de Physique Th\'{e}orique, University 
of Geneva, Spring 2001, during which some of this work was done. He was
supported in part by NSF Grant 980139. A.B. was supported in part by
NSF grant 9971493 and by a Sloan fellowship.


\end{document}